\documentclass[review,hidelinks,onefignum,onetabnum]{siamart251216}

\usepackage{lipsum}
\usepackage{amsfonts}
\usepackage{graphicx}
\usepackage{epstopdf}
\usepackage{algorithmic}

\usepackage{graphicx,subfigure}
\usepackage{epstopdf}
\usepackage{nicefrac}
\usepackage{bm}
\usepackage{bbm}
\usepackage{amssymb}
\usepackage{amsmath}
\usepackage{mathrsfs}
\usepackage{mathtools}
\usepackage{color}
\usepackage{stmaryrd}
\usepackage[colorinlistoftodos,prependcaption,textsize=tiny]{todonotes}
\usepackage{extarrows}
\usepackage{xspace}
\usepackage[normalem]{ulem} 

\ifpdf
  \DeclareGraphicsExtensions{.eps,.pdf,.png,.jpg}
\else
  \DeclareGraphicsExtensions{.eps}
\fi

\newcommand{\dt}{\Delta t}
\newcommand{\mba}{\mathbf{a}}
\newcommand{\mbx}{\mathbf{x}}
\newcommand{\mbu}{\mathbf{u}}
\newcommand{\mbv}{\mathbf{v}}
\newcommand{\mbf}{\mathbf{f}}
\newcommand{\mbJ}{\mathbf{J}}
\newcommand{\mbn}{\mathbf{n}}

\newcommand{\Id}{\mathrm{Id}}

\newcommand{\R}{\mathbb{R}}
\newcommand{\V}{\mathbb{V}}
\newcommand{\N}{\mathcal{N}}

\newcommand{\dpar}[2]{\dfrac{\partial #1}{\partial #2}}

\newcommand*\xbar[1]{%
  \hbox{%
    \vbox{%
      \hrule height 0.5pt 
      \kern0.4ex
      \hbox{%
        \kern-0.05em
        \ensuremath{#1}%
        \kern-0.00em
      }%
    }%
  }%
}

\newcommand{\DD}{\mathcal D}

\newcommand{\II}{\mathcal I}

\newcommand{\PP}{\mathcal P}

\newcommand{\bba}{\mathbf{a}}

\newcommand{\bbf}{\mathbf{f}}

\newcommand{\bbn}{\mathbf{n}}

\newcommand{\bbu}{\mathbf{u}}
\newcommand{\bbv}{\mathbf{v}}

\newcommand{\bbx}{\mathbf{x}}

\newcommand{\dg}{DG\xspace}
\newcommand{\pampa}{PAMPA \xspace}
\newcommand{\red}[1]{\textcolor{red}{#1}}
\newcommand{\blue}[1]{\textcolor{blue}{#1}}

\hypersetup{
    colorlinks=true,
    linkcolor=green,     
    citecolor=red,      
    urlcolor=blue    
}

\newsiamremark{remark}{Remark}
\newsiamremark{hypothesis}{Hypothesis}
\crefname{hypothesis}{Hypothesis}{Hypotheses}
\newsiamthm{claim}{Claim}
\newsiamremark{fact}{Fact}
\crefname{fact}{Fact}{Facts}

\headers{Robust \pampa Scheme in the DG Formulation}{R. Abgrall and Y. Liu}

\graphicspath{{figures/}{./}}

\title{Robust \pampa Scheme in the DG Formulation on Unstructured Triangular Meshes: Bound Preservation, Oscillation Elimination, and Boundary Conditions \thanks{Submitted to the editors DATE.
\funding{The work of Y. Liu was supported by UZH Postdoc Grant, 2024 / Verf\"{u}gung Nr. FK-24-110 and SNSF grant 200020$\_$204917.}}}

\author{R\'emi Abgrall\thanks{Institute of Mathematics, University of Z\"{u}rich, 8057 Z\"{u}rich, Switzerland (\email{remi.abgrall@math.uzh.ch})} \and Yongle Liu\thanks{Corresponding Author. Institute of Mathematics, University of Z\"{u}rich, 8057 Z\"{u}rich, Switzerland
  (\email{yongle.liu@math.uzh.ch}).}}

\usepackage{amsopn}

\begin{document}
\nolinenumbers
\maketitle

\begin{abstract}
We propose an improved version of the PAMPA algorithm where the solution is sought as globally continuous. The scheme is locally conservative, and there is no mass matrix to invert. This method had been developed in a series of papers, see e.g \cite{Abgrall2024a} and the references therein. In \cite{Abgrall2025d}, we had shown the connection between \pampa and the discontinuous Galerkin method, for the linear hyperbolic problem. Taking advantage of this reinterpretation, we use it to define a family of methods, show how to implement the boundary conditions in a rigorous manner. In addition, we propose a method that complements the bound preserving method developed in \cite{Abgrall2025d} in the sense that it is non oscillatory. A truncation error analysis is provided, it shows that the scheme should be third order accurate for smooth solutions. This is confirmed by numerical experiments. Several numerical examples are presented to show that the scheme is indeed  bound preserving and non oscillatory on a wide range on numerical benchmarks.
\end{abstract}



\section{Introduction}\label{sec1}
Since the seminal works of P. L. Roe and his students, \cite{Eymann2011,Eymann2011a,Eymann2013,Eymann2013a,Maeng2017,Fan2017,He2021}, there has been a growing interest in the so-called Active Flux (AF) schemes for solving hyperbolic problems. Unlike Godunov-type methods, the AF schemes evolve simultaneously the average values and point values located at the boundaries of elements. The cell averages are updated through a conservative formulation with the numerical fluxes being computed directly from the point values. The evolution of point values is more involved and two main strategies have been proposed. The first employs exact or approximate evolution operators \cite{Chudzik2024,Chudzik2021,Barsukow2024,Helzel2019,Barsukow2021}. However, it is often difficult to construct efficient approximated operators for multi-dimensional or nonlinear systems without degrading high order accuracy. The second approach, initiated by Abgrall \cite{Abgrall2023a} and known as the semi-discrete AF method, avoids the need for such operators by first using the degrees of freedom (DoFs) to discretize the spatial derivatives, and then applying the method-of-line to integrate the resulting semi-discrete system. Much progress has been made in the development and application of the semi-discrete AF method. Notable contributions include its arbitrary high-order extensions in one-dimension (1D) settings \cite{Abgrall2023b,Abgrall2023c}, applications to Euler equations on Cartesian grids \cite{Abgrall2023}, compressible flows on unstructured triangular meshes \cite{Abgrall2025a} and arbitrary polygonal meshes \cite{Abgrall2024a}, and the enforcement of bound-preserving (BP) properties in both 1D \cite{Abgrall2026,Abgrall2025b} and two-dimensional (2D) \cite{Abgrall2024a,Duan2025} cases. Moreover, this methodology has also been extended to 1D hyperbolic balance laws \cite{Abgrall2024,Liu2025} and 2D hyperbolic balance laws \cite{Liu2025a}, where it is referred to as the \pampa (Point-Average-Moment PolynomiAl-interpreted) scheme. Several very intriguing properties of the \pampa (semi-discrete AF) have been noticed. In contrast to the traditional AF method, the computational elements in \pampa can be intervals (not necessarily uniform) in 1D, and quads, triangles, or even general polygons in several space dimensions. Furthermore, the evolution of point values is not restricted to conservative formulations; instead, it allows customized updates based on arbitrary non-conservative formulations tailored to specific requirements. For example, the primitive variables used in \cite{Abgrall2023a,Abgrall2024,Liu2025}, the equilibrium variables adopted in \cite{Liu2025a}, and the unconditionally limiter-free variables used in \cite{Abgrall2025b}. This added flexibility makes \pampa particularly suitable for complex hyperbolic systems. In addition, it seems to have very interesting properties for the low Mach flows \cite{Abgrall2024a}, which, in our opinion, has been partially explained, at least in the Cartesian version \cite{Barsukow2025}. Finally, as investigated in \cite{Abgrall2025d}, by choosing an appropriate projection step from the discontinuous approximation space back to the continuous one, a connection can be established between \pampa and the Discontinuous Galerkin (\dg) methods.

This paper is a follow-up of \cite{Abgrall2025d}, where we have shown that the 1D scheme developed in \cite{Abgrall2023a,Abgrall2026} can be reinterpreted as a combination of the \dg method together with a projection step on the family of globally continuous piecewise polynomials. In 1D, several projectors have already been identified. The one that is implicitly used in \cite{Abgrall2023a,Abgrall2026} employs an upwinding mechanism, while the alternative presented in \cite{Abgrall2025d} relies on an arithmetic averaging approach. Unfortunately, numerical experiments in 2D, particularly on unstructured triangular meshes, indicate that the latter (arithmetic average) projector leads to instabilities after a sufficiently long integration time, whereas the upwind-type projector remains stable, at least experimentally. This last one is very close to what has already been proposed in \cite{Abgrall2025a,Abgrall2024a}. It is suggested that a large family of projectors can make the job but using a combination of upwinding and centering, a bit as in the spirit of upwind biased finite differences. 

In this  paper, we further demonstrate the robustness of the \pampa scheme, based on a \dg formulation combined with an upwind-type projector, through applications to scalar conservation laws and the compressible Euler equations. 
We also present a consistent discretization of boundary conditions, made possible by the variational interpretation of the method. For the sake of completeness, independently and at the same time, the submission \cite{connard} proposed a similar interpretation of the AF method on Cartesian grids. Furthermore, We present a truncation error analysis in Appendix \ref{sec:Truncation_error}, explaining why the scheme is third order accurate. Note that this analysis can be easily generalized to higher order approximations obtained following \cite{Abgrall2023b,Abgrall2023c}.

We are interested in hyperbolic problems of the type:
\begin{equation}
\label{eq:1}
\dpar{\mbu}{t}+\text{ div }\mbf(\mbu)=0, \quad t>0, \quad \mbx\in \Omega\subset \R^d
\end{equation}
where $\mbu\in \DD\subset \R^m$, and $\mbf=(f_1, \ldots , f_d)$ is assumed to be $C^1$. The set $\Omega$ is open in $\R^d$. The domain $\DD$, called the invariant domain, is assumed to be convex and open in $\R^m$. 
Here, \eqref{eq:1} is complemented by an initial condition 
\begin{equation*}
  \mbu(\mbx,0)=\mbu_0(\mbx)
\end{equation*}
and boundary conditions that will be described in the relevant section. Without loss of generality, we focus on the following 2D hyperbolic equations:
\begin{itemize}
\item the scalar conservation laws with  linear flux $\mbf$ of  the form  $\mbf(\mbu,\mbx)=\mba(\mbx)\mbu$ or a non-linear flux. In that case, $\mbu$ is a function with values in $\R$ and $\DD=\R$ but the solution must stay in $[\min\limits_{\mbx\in \R}\mbu_0(\mbx), \max\limits_{\mbx\in \R}\mbu_0(\mbx)]$, due to Kruzhkov's theory. In these particular examples, we set $\DD=[\min\limits_{\mbx\in \R}\mbu_0(\mbx), \max\limits_{\mbx\in \R}\mbu_0(\mbx)]$.
\item the Euler equations with the perfect gas equation of state, and
\begin{equation*}
     \mbu=\begin{pmatrix}
         \rho \\
          \rho \mbv\\
           E
      \end{pmatrix}, \quad 
\mbf(\mbu)=\begin{pmatrix}
\rho \mbv\\
\rho \mbv\otimes \mbv+p\text{Id}_d\\
(E+p)\mbv
\end{pmatrix}.
\end{equation*}
Here, $\text{Id}_d$ is a $d\times d$ identity matrix, $\rho$ is the density, $\mbv=(u,v)$ is the velocity field, and $E$ is the total energy, i.e., the sum of the internal energy $\epsilon$ and the kinetic energy $\tfrac{1}{2}\rho \mbv^2$. We have introduced the pressure that is related to the internal energy and the density by an equation of state. Here we make the choice of a perfect gas, thus
\begin{equation*}
 p=(\gamma-1)\epsilon
\end{equation*}
with $\gamma$ denoting the specific heat ratio. In this example, the invariant domain $\DD=\{\mbu\in \R^4, \rho>0, \epsilon>0\}$. Thanks to the geometric quasi-linearization (GQL) representation given in \cite{Wu2023}, the equivalent form of the invariant domain $\mathcal{D}$ is given by
\begin{equation}\label{eq:2}
   \mathcal{D}_{\bm \nu}=\{\mbu=(\rho, \rho\mbv, E)^\top\in\R^4 \text{ such that for all } \mbn_*\in \N, \mbu^\top\mbn_*>0\}, 
\end{equation}
where the vector space $\N$ is given by
\begin{equation*}
     \mathcal{N}=\left\{\begin{pmatrix} 1\\ \mathbf{0}_d\\0\end{pmatrix}, \mathbf{0}_d\in\R^d\right\}\cup\left\{ \begin{pmatrix}\frac{\Vert\bm\nu\Vert^2}{2} \\ -\bm\nu\\ 1\end{pmatrix}, \bm\nu\in \R^d\right\}.
\end{equation*}
\end{itemize}

It is known that, high-order numerical methods for nonlinear hyperbolic PDEs are prone to generating unwanted spurious oscillations near shocks and discontinuities. These oscillations can result in nonphysical values such as negative density and/or pressure, leading to code break down. 

To solve this issue, we have derived a provably BP \pampa scheme in \cite{Abgrall2026,Abgrall2024a} based on the monolithic convex limiting approach. In that case, the density and the internal energy will numericaly stay above a threshold $\epsilon>0$ that can be arbitrarily chosen. Conservation is not broken.  It was further shown that, by employing the equivalent invariant domain \eqref{eq:2}, the explicit BP coefficients can be obtained from a simple quadratic formulation. Therefore, in this paper, to ensure the BP property, we directly employ the BP parameters derived in \cite{Abgrall2024a,Abgrall2026} to blend the high- and low-order fluxes and residuals for updating internal and boundary DoFs. \red{Other techniques to ensure the invariant domain preserving, we refer a recent review paper \cite{Wu2025} and the references therein.}

However, this is not enough: a BP limiter will not prevent spurious oscillations, it will only prevent that the density and/or the internal energy will not become negative.

To mitigate these oscillations, a variety of nonlinear stabilization techniques, often referred to as shock capturing methods, have been developed. They can be broadly categorized into one of four groups: limiting \cite{Cockburn1990,Cockburn1989,Krivodonova2007,Krivodonova2004,Lu2019,Hoteit2004,Park2012,Park2014,Park2016,VanLeer1977,Woodward1984,Yang2009,You2018}, filtering \cite{Asthana2015,Hesthaven2008}, artificial viscosity \cite{Fu2022,Kloeckner2011,Park2014a,Persson2006,Yang2016,Huang2020}, or method modification \cite{Dumbser2014,Loubere2014,Boscheri2015,Diot2013}. 
Recently, a novel approach to control spurious oscillations by introducing damping terms into the semi-discrete \dg scheme, termed the oscillation-free \dg (OFDG, \cite{Lu2021}) and oscillation-eliminating \dg (OEDG, \cite{Peng2025}) methods have been successfully applied to various hyperbolic problems (e.g., \cite{Ding2025,Fan2024,Peng2025a,Yan2024,Liu2025c,Liu2024}). In particular, the OEDG method features a novel scale- and evolution-invariant damping operator that avoids characteristic decomposition and effectively eliminates spurious oscillations in problems with varying scales and wave speeds. Very recently, a new convex oscillation-suppressing (COS) framework for high-order \dg methods is proposed \cite{Cao2026}. The idea of COSDG is to integrate the oscillation elimination parameter and bound preservation parameter in a unified framework. Inspired by these filtering techniques, especially the OEDG and COSDG methods, we incorporate a damping term based on high-order derivative jumps at element interfaces into the 1D BP \pampa scheme. This leads to the development of a non-intrusive convex limiting parameter designed to suppress spurious oscillations, as introduced in \cite{Abgrall2025c}. To elaborate, we define a convex oscillation eliminating parameter in $[0,1]$ using the damping term. After combining it with the BP parameters via minimization, this blended parameter is used to perform convex blending between the high- and low-order fluxes and residuals. 

The contribution of this paper is to present the BP OE \pampa scheme in such a way that the \dg method is one of the building block. In order to control spurious oscillations, our recent 1D convex oscillation eliminating procedure \cite{Abgrall2025c} (inspired by \cite{Abgrall2025b,Peng2025,Lu2021,Chudzik2025}) has been extended here on unstructured triangular meshes.
The main findings of this work can be summarized as follows:
\begin{itemize}
       \item We present a rigorous derivation of the \pampa scheme within a \dg framework. The key steps involve constructing the mass matrix using the same polynomial basis functions, applying the Riesz representation theorem to express linear forms as inner products between primal and dual basis functions, and formulating the dual basis within the same finite element polynomial space. To recover the \pampa scheme from its \dg formulation, we introduce a suitable projection operator. While the evolution of internal DoF is retained in its original \dg form, the evolution of boundary DoFs must be modified due to the global continuity assumption inherent in the standard \pampa framework. Drawing inspiration from Residual Distribution schemes, we express the update of boundary DoFs as a convex combination with appropriately chosen weights. In particular, we provide a choice of weights consistent with the upwind \pampa approach introduced in \cite{Abgrall2025a,Abgrall2024a}. 

            \item To guarantee the BP property of the \pampa scheme, we adopt the monolithic convex limiting strategy introduced in \cite{Abgrall2026,Abgrall2024a}, where the BP parameter is explicitly computed using GQL representation \eqref{eq:2} for the Euler equations. To effectively suppress spurious oscillations near strong discontinuities and shock waves, building upon prior work for 1D shallow water flows in \cite{Abgrall2025c}, we generalize the convex oscillation eliminating parameter to multidimensional. The key idea is to measure discrepancies in numerical solutions across neighboring elements, thereby identifying regions exhibiting oscillatory behavior. This measurement is taken from the OFDG and OEDG methods, which utilize jumps of higher-order derivatives across shared element interfaces. 
                In the end, the OE parameter is then combined with the BP parameter by taking their minimum, resulting in a convex BP OE limiting parameter. This parameter is used to blend the high- and low-order numerical fluxes and residuals in a non-intrusive and robust fashion.
                
         \item It is worth noting that the treatment of boundary conditions within the \pampa scheme is nontrivial, particularly on unstructured meshes \cite{Abgrall2024a}. The \dg formulation offers a systematic and consistent framework for implementing boundary conditions. In this approach, the boundary contributions naturally can be incorporated directly into the residuals used to update both internal and boundary DoFs. This greatly simplifies the implementation of the scheme for various boundary conditions.
\end{itemize}

{\color{blue} Finally, one is legitimate to ask why a new scheme and why not simply doing \dg? One answer is to count the number of degrees of freedom in 3D. From the Euler relations, a regular mesh made of tetrahedrons with $n_v$ vertices has about $7n_v$ edges, $12 n_v$ and $6n_v$ tetrahedrons. Then a quick calculation allows to evaluate the number of degrees of freedom in the case of \dg and a globally continuous representation of the solution, see Figure \ref{ratio}. This shows that for moderate order, the continuous finite element methods are more competitive in term of storage compared with \dg and the type of method we are advocating, to the price of a global mass matrix to handle. In the class of method without global mass matrix, then we see that the method we are advocating are more competitive than \dg. Of course this is not the only criteria, but this shows there is room for research. 
\begin{figure}[h]
\begin{center}
\subfigure[]{\includegraphics[width=0.45\textwidth]{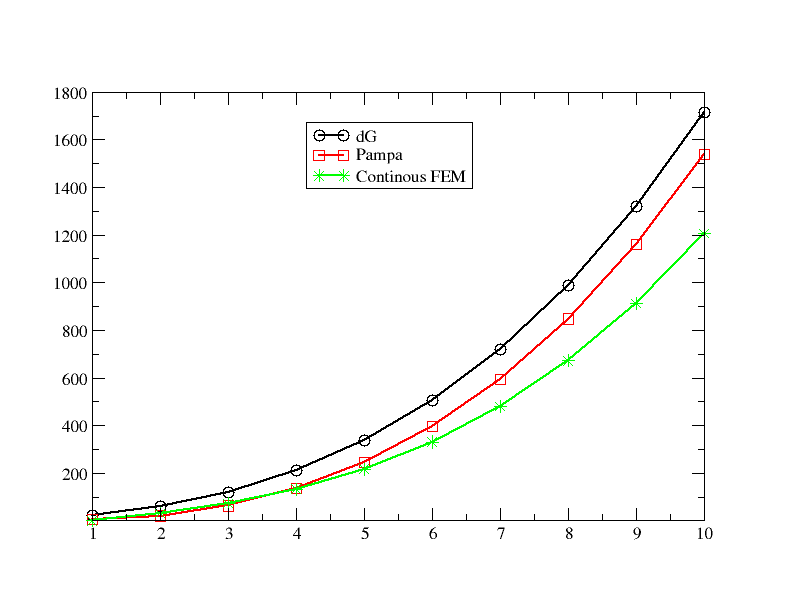}}
\subfigure[]{\includegraphics[width=0.45\textwidth]{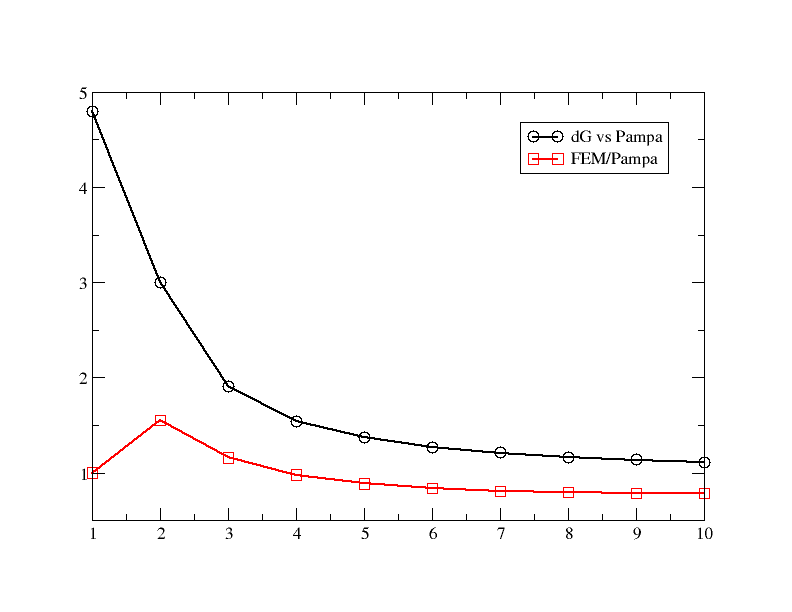}}
\caption{\label{ratio}For 3D tetrahedron meshes, (a): Ratio of the number of degrees of freedom with the number of vertices, (b) Ration of the number of degrees of freedom between \dg, continuous FEM and \pampa-like methods.}
\end{center}
\end{figure}
In two dimensions for triangular elements and  order equal or larger than $3$, by extending  the work of \cite{Abgrall2025a}, it is possible to have the same number of degrees of freedoms as for a continuous finite element method of the same order. It is possible that the same conclusion holds also in 3D, but this needs to be investigated. The extension to non triangular/tetrahedron elements also need to be investigated.}

The format of this paper is as follows. First, we recall the relevant elements of \pampa schemes in \cite{Abgrall2025d}, and establish the connection between the \pampa and \dg methods through the Riesz representation theorem, as well as define the projector operators so that the upwind \pampa scheme introduced in \cite{Abgrall2025a,Abgrall2024a} can be recovered from its \dg formulation. Then we provide a brief review of the low-order \pampa schemes that can be used to ensure BP and OE properties. Next, we introduce the BP parameters in \cite{Abgrall2024a} and propose a rotational invariant convex oscillation eliminating parameter using the damping term computed by high-order derivatives of the approximated solutions. Last, we discuss the boundary condition problem and show several challenging numerical examples on the scalar case and compressible Euler equations of gas dynamics. 

\section{Robust \pampa Schemes}\label{sec2}
In this section, we review the essential ingredients of the standard high-order \pampa scheme on 2D unstructured triangular meshes, $\Omega=\bigcup_j K_j$.We then demonstrate how the main building block of the high-order \pampa scheme can be interpreted as a conventional \dg method, acting on globally continuous functions. Next, we describe the low-order \pampa scheme constructed in \cite{Abgrall2025a}, as well as the BP parameters introduced in \cite{Abgrall2024a}. Finally, in order to control spurious oscillations near strong discontinuities, we propose an element-based oscillation-eliminating (OE) parameter, which is inspired by \cite{Abgrall2025b,Peng2025,Lu2021} and defined using the jumps of derivatives of the approximate numerical solutions. The OE parameter, combined with the BP parameter, will be used to define a convex blending of the high-order and low-order \pampa schemes, ensuring that the resulting \pampa scheme is both BP and OE robust.

\begin{figure}[ht!]\label{fig:1}
\centerline{\includegraphics[trim=0.01cm 0.01cm 0.01cm 0.001cm,clip,width=4.0cm]{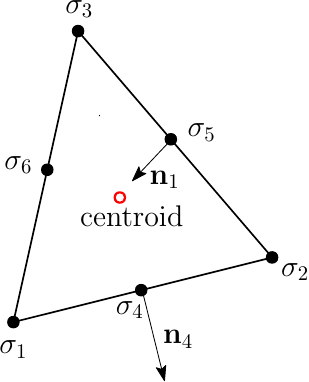}}
\caption{DoFs and normals of triangle $K$.}
\end{figure}

For the sake of simplicity, we reduce ourselves to the third-order case, where DoFs consist of boundary DoFs---6 point values located at the three vertices and the three edge midpoints, and one internal DoF corresponding to the cell average within each element; see Figure \ref{fig:1}. In each triangle $K$, the vertices are denoted by $\sigma_1$, $\sigma_2$, and $\sigma_3$, and the barycentric coordinates by $(\lambda_1, \lambda_2, \lambda_3)$. The other degrees of freedom are the mid-points of the edges, denoted by $\sigma_4,\sigma_5,\sigma_6$. For higher-order extensions, several strategies to increase the number of DoFs have been proposed in \cite{Abgrall2023b} for 1D intervals. Among these, we favor the approach with additional higher-order moments, as it preserves the compactness of the scheme and facilitates its generalization to multiple spatial dimensions. The underlying idea for both 1D and multi-D settings has been discussed in \cite{Abgrall2025d}, while the detailed reconstruction procedure and related applications are left for future work.

\subsection{Third-order \pampa Scheme and Its \dg Formulation}
In the standard third-order \pampa (or the so-called semi-discrete/generalized AF) scheme, the solution of \eqref{eq:1} is approximated by a globally continuous finite element polynomial expansion $\mbu_{\mathrm h}$ within each element $K$:
\begin{equation}\label{eq:3}
  \mbu_{\mathrm h}|_K=\sum_{i=1}^{6}\mbu_{\sigma_i}\varphi_{\sigma_i}+\xbar{\mbu}_K\xbar\varphi,
\end{equation}
where the basis functions are given by (\cite{Abgrall2025a,Abgrall2025d})
\begin{itemize}
\item For the average, we set $\xbar\varphi=60\lambda_1\lambda_2\lambda_3$,
\item For the vertices $\{\sigma_i\}_{i=1}^3$, we set
\begin{equation*}
  \varphi_i=\lambda_{i}\big (2\lambda_{i}-1),
\end{equation*}
\item and for the midpoint of the edge $[\sigma_i,\sigma_{i+1}]$ modulo 3, denoted by $\sigma_{3+i}$, we set
\begin{equation*}
  \varphi_{3+i}=4\lambda_{i}\lambda_{{i+1}}-\frac{1}{3}\xbar\varphi.
\end{equation*}
\end{itemize}
We, therefore, construct the local approximation space $\V_k(K)=\text{span}\{\varphi_i,\xbar\varphi\}$ for each triangle element $K$. Following \cite{Abgrall2025d}, we first introduce the linear forms ($\V_k(K)\rightarrow\R$):  
\begin{equation*}
  \text{DoF}_i(\mbu_{\mathrm h})=i\text{-th DoF of } \mbu_{\mathrm h}, \quad i=1,\ldots, N_{\text{DoFs}}.
\end{equation*}
Here, $N_{\text{DoFs}}=\text{dim} \V_k(K)=7$. We can then write \eqref{eq:3} in terms of these linear forms by means of
\begin{equation}\label{eq:uh}
  \mbu_{\mathrm h}|_K=\sum_{i=1}^{N_{\text{DoFs}}}\text{DoF}_i(\mbu_{\mathrm h})\varphi_i,
\end{equation}
where we have specified $\varphi_7=\xbar \varphi$ and the cell average $\xbar\mbu_K$ is the corresponding $7$-th DoF, i.e.,
\begin{equation*}
  \text{DoF}_7(\mbu_{\mathrm h})=\xbar\mbu_K=\frac{1}{\vert K\vert}\int_K \mbu_{\mathrm h}(\mbx)\; \mathrm d\mbx.
\end{equation*}
Furthermore, we have the usual interpolation property holds true:
\begin{equation}
  \text{DoF}_i(\varphi_j)=\delta_{ij},\quad i,j=1,\ldots,N_{\text{DoFs}}.
\end{equation} 
Using Riesz theorem with  the scalar product between $\bbu$ and $\bbv\in \V_k(K)$ defined by
$$\langle \bbu,\bbv\rangle=\frac{1}{\vert K\vert }\int_K \bbu(\bbx)\cdot \bbv(\bbx)\; \mathrm d\bbx, $$ we can consider the dual basis $\theta_{\sigma_i}\in\V_k(K)$ for $i=1, \ldots, 6$, such that
\begin{equation}\label{eq:DoF_sigma}
  \text{DoF}_i(\mbu_{\mathrm h})=\mbu_{\mathrm h}(\sigma_i)=\frac{1}{\vert K\vert}\int_K \theta_{\sigma_i}(\mbx)\cdot\mbu_{\mathrm h}(\mbx)\; \mathrm d\mbx.
\end{equation}
and for the $7$-th DoF (cell average) we should obviously choose the dual basis $\theta_{\mu}\in\V_k(K)$ as $1$ to obtain 
\begin{equation}\label{eq:DoF_mu}
  \text{DoF}_7(\mbu_{\mathrm h})=\xbar\mbu_K=\frac{1}{\vert K\vert}\int_K \theta_{\mu}(\mbx)\cdot\mbu_{\mathrm h}(\mbx)\; \mathrm d\mbx=\frac{1}{\vert K\vert}\int_K \mbu_{\mathrm h}(\mbx)\; \mathrm d\mbx.
\end{equation}
By writing 
\begin{equation*}
  \theta_{\sigma_i}=\sum_{l=1}^7 a_{\sigma_i,l}\varphi_l,\quad  \theta_{\mu}=\sum_{l=1}^7 a_{\mu,l}\varphi_l,
\end{equation*}
and inserting into
\begin{equation*}
\begin{aligned}
  &\text{DoF}_i(\varphi_j)=\delta_{ij}=\frac{1}{\vert K\vert}\int_{K} \theta_{\sigma_i}\varphi_j\;\mathrm d\mbx,\quad i=1,\ldots, 6,~j=1,\ldots,7\\
  &\text{DoF}_7(\varphi_j)=\delta_{7j}=\frac{1}{\vert K\vert}\int_{K} \theta_{\mu}\varphi_j\;\mathrm d\mbx,\quad j=1,\ldots,7,
\end{aligned}
\end{equation*}
we see that
\begin{equation}\label{eq:dual_basis_ex}
\begin{pmatrix}
 \theta_{\sigma_1}\\
 \vdots\\
 \theta_{\sigma_6}\\
 \theta_\mu
\end{pmatrix}=\vert K\vert M_K^{-1}\begin{pmatrix}
\varphi_1 \\ \vdots \\
\varphi_6\\
\varphi_7
\end{pmatrix},
\end{equation}
where $[M_K]_{ij}=\int_K \varphi_i\varphi_j\,\mathrm d\mbx$ is the mass matrix. A detailed proof, for the general case, we refer the readers to \cite{Abgrall2025d}.
In our case (third-order and triangular meshes), we get
\begin{equation}\label{eq:Mk_inv}
  M_K^{-1}=\frac{1}{\vert K\vert}\begin{pmatrix}
  &\frac{140}{3} &\frac{50}{3}&\frac{50}{3}&-\frac{65}{6}&-\frac{10}{3}&-\frac{65}{6}&1\\[0.3em]
  &\frac{50}{3} &\frac{140}{3}&\frac{50}{3}&-\frac{65}{6}&-\frac{65}{6}&-\frac{10}{3}&1\\[0.3em]
  &\frac{50}{3} &\frac{50}{3}&\frac{140}{3}&-\frac{10}{3}&-\frac{65}{6}&-\frac{65}{6}&1\\[0.3em]
  &-\frac{65}{6} &-\frac{65}{6}&-\frac{10}{3}&\frac{215}{12}&\frac{115}{24}&\frac{115}{24}&1\\[0.3em]
  &-\frac{10}{3} &-\frac{65}{6}&-\frac{65}{6}&\frac{115}{24}&\frac{215}{12}&\frac{115}{24}&1\\[0.3em]
  &-\frac{65}{6} &-\frac{10}{3} &-\frac{65}{6}&\frac{115}{24}&\frac{115}{24}&\frac{215}{12}&1\\[0.3em]
  &1&1&1&1&1&1&1
  \end{pmatrix},
\end{equation}
and thus
\begin{equation*}
\begin{split}
\theta_{\sigma_1}&=\frac{140 \lambda_1 \left(2 \lambda_1 -1\right)}{3}+\frac{50 \lambda_2 \left(2 \lambda_2 -1\right)}{3}+\frac{50 \lambda_3 \left(2 \lambda_3 -1\right)}{3}+560 \lambda_1 \lambda_2 \lambda_3 -
\\
&\qquad\frac{130 \lambda_1 \lambda_2}{3}-\frac{40 \lambda_2 \lambda_3}{3}-\frac{130 \lambda_3 \lambda_1}{3}, 
\\
 \theta_{\sigma_2}&=\frac{50 \lambda_1 \left(2 \lambda_1 -1\right)}{3}+\frac{140 \lambda_2 \left(2 \lambda_2 -1\right)}{3}+\frac{50 \lambda_3 \left(2 \lambda_3 -1\right)}{3}+560 \lambda_1 \lambda_2 \lambda_3 -
\\
&\qquad\frac{130 \lambda_1 \lambda_2}{3}-\frac{130 \lambda_2 \lambda_3}{3}-\frac{40 \lambda_3 \lambda_1}{3}, 
\end{split}
\end{equation*}
\begin{equation*}
    \begin{split}
\theta_{\sigma_3}&=\frac{50 \lambda_1 \left(2 \lambda_1 -1\right)}{3}+\frac{50 \lambda_2 \left(2 \lambda_2 -1\right)}{3}+\frac{140 \lambda_3 \left(2 \lambda_3 -1\right)}{3}+560 \lambda_1 \lambda_2 \lambda_3 -
\\
&\qquad\frac{40 \lambda_1 \lambda_2}{3}-\frac{130 \lambda_2 \lambda_3}{3}-\frac{130 \lambda_3 \lambda_1}{3}, 
\\
\theta_{\sigma_4}&= -\frac{65 \lambda_1 \left(2 \lambda_1 -1\right)}{6}-\frac{65 \lambda_2 \left(2 \lambda_2 -1\right)}{6}-\frac{10 \lambda_3 \left(2 \lambda_3 -1\right)}{3}-490 \lambda_1 \lambda_2 \lambda_3 +
\\
&\qquad\frac{215 \lambda_1 \lambda_2}{3}+\frac{115 \lambda_2 \lambda_3}{6}+\frac{115 \lambda_3 \lambda_1}{6}, 
\\
\theta_{\sigma_5}&= -\frac{10 \lambda_1 \left(2 \lambda_1 -1\right)}{3}-\frac{65 \lambda_2 \left(2 \lambda_2 -1\right)}{6}-\frac{65 \lambda_3 \left(2 \lambda_3 -1\right)}{6}-490 \lambda_1 \lambda_2 \lambda_3 +
\\
&\qquad\frac{115 \lambda_1 \lambda_2}{6}+\frac{215 \lambda_2 \lambda_3}{3}+\frac{115 \lambda_3 \lambda_1}{6}, 
\\
\theta_{\sigma_6}&= -\frac{65 \lambda_1 \left(2 \lambda_1 -1\right)}{6}-\frac{10 \lambda_2 \left(2 \lambda_2 -1\right)}{3}-\frac{65 \lambda_3 \left(2 \lambda_3 -1\right)}{6}-490 \lambda_1 \lambda_2 \lambda_3 +
\\
&\qquad\frac{115 \lambda_1 \lambda_2}{6}+\frac{115 \lambda_2 \lambda_3}{6}+\frac{215 \lambda_3 \lambda_1}{3}, 
\\
\theta_{\mu}&=1.
\end{split}
\end{equation*}
We note that, from \eqref{eq:Mk_inv},
\begin{equation}\label{eq:Project}
  M_K^{-1}=P_K=\frac{1}{|K|} \mathcal{P},
\end{equation}
where the matrix $\mathcal{P}$ does not depend on $K$, because we are using barycentric coordinate: this is as if working in the reference element. 

On the other hand, using the same basis $\{\varphi_i\}_{i=1}^7$, the \dg scheme in one element $K$ writes
\begin{equation}\label{eq:DG_form}
  \dfrac{\mathrm d}{\mathrm dt}\mathbf U_K+M_K^{-1}\mathbf F_K=0,
\end{equation}
where 
\begin{equation}\label{eq:dg_vectors}
  \mathbf U_K=\begin{pmatrix}
          \mathbf U_{K,\sigma}\\
          \mathbf U_{K,\mu}
  \end{pmatrix}=\begin{pmatrix}
                \mbu_{\sigma_1} \\
                \mbu_{\sigma_2} \\
                \vdots \\
                 \xbar\mbu_K
              \end{pmatrix},\quad
  \mathbf F_K=\begin{pmatrix}
          \mathbf F_{K,\sigma}\\
          \mathbf F_{K,\mu}
  \end{pmatrix}=\begin{pmatrix}
                F_{\sigma_1} \\
                F_{\sigma_2} \\
                \vdots \\
                 F_\mu
              \end{pmatrix}
\end{equation}
with 
\begin{equation}\label{eq:dg_load}
\begin{aligned}
F_{\sigma_i}&=
    -\int_{K} \nabla \varphi_{i} \cdot \mbf(\mbu_{\mathrm h})\; \mathrm d\mbx+\int_{\partial K}\varphi_{i}\mbf(\mbu_{\mathrm h})\cdot \mbn \; \mathrm d\ell,\quad i=1,\ldots,6,\\
F_\mu&=-\int_{K} \nabla \varphi_{7} \cdot \mbf(\mbu_{\mathrm h})\; \mathrm d\mbx+\int_{\partial K}\varphi_{7}\mbf(\mbu_{\mathrm h})\cdot \mbn \; \mathrm d\ell.
\end{aligned}
\end{equation}
To interpret the \pampa scheme in the \dg formulation \eqref{eq:DG_form}--\eqref{eq:dg_load}, let us apply the linear forms in \eqref{eq:DoF_sigma} and \eqref{eq:DoF_mu} onto the PDE \eqref{eq:1} and obtain
\begin{equation}\label{eq:PAMPA_DG}
\begin{aligned}
  &\frac{\mathrm d}{\mathrm dt}\int_K\theta_{\sigma_i}\cdot\mbu_{\mathrm h}(\mbx)\;\mathrm d\mbx+\int_K\theta_{\sigma_i}\cdot{\text{div}}~\mbf\;\mathrm d\mbx=0,\quad i=1,\ldots,6\\
  &\frac{\mathrm d}{\mathrm dt}\int_K\theta_{\mu}\cdot\mbu_{\mathrm h}(\mbx)\;\mathrm d\mbx+\int_K\theta_{\mu}\cdot{\text{div}}~\mbf\;\mathrm d\mbx=0.
  \end{aligned}
\end{equation}
This yields
\begin{equation*}
  \dfrac{\mathrm d}{\mathrm dt}\mathbf U_K+\frac{1}{\vert K\vert}
  \int_K \begin{pmatrix}
  \{\theta_{\sigma_i}\}_{i=1}^6 \\
  \theta_\mu
  \end{pmatrix}\cdot{\text{div}}~\mbf\;\mathrm d\mbx=0.
\end{equation*}
Combing with \eqref{eq:dual_basis_ex}, we have
\begin{equation}\label{equation:dg}
  \dfrac{\mathrm d}{\mathrm dt}\mathbf U_K+
  M_K^{-1}\int_K \begin{pmatrix}
  \{\varphi_{i}\}_{i=1}^6 \\
  \varphi_7
  \end{pmatrix}\cdot{\text{div}}~\mbf\;\mathrm d\mbx=0.
\end{equation}
Because everything is assumed to be globally smooth in $K$ and the divergence theorem, we exactly get the \dg formulation \eqref{eq:DG_form}--\eqref{eq:dg_load}.
Let us remark that the evolution equation of the average $\xbar\mbu_K$ coming out of \eqref{equation:dg} is simply
\begin{equation}\label{equation:average}
\dfrac{\mathrm d\xbar \mbu_K}{\mathrm dt}+\frac{1}{\vert K\vert}\int_{\partial K} \mbf(\mbu_{\mathrm h})\cdot \mbn\; \mathrm d\ell=0.
\end{equation}

\subsection{From \dg to \pampa}
Here, we discuss how to recover the update of high-order \pampa scheme from its \dg formulation \eqref{eq:DG_form}--\eqref{eq:dg_load}. We partition the time domain in intermediate times $\{t^n\}$ with $\dt^n=t^{n+1}-t^n$ the $n$-th time step and denote the solution at time $t^n$ by $\mbu^n$ and for the sake of simplicity, we consider the forward Euler time stepping from $t^n$ to $t^{n+1}$ for solving \eqref{eq:DG_form}, which reads as
\begin{subequations}\label{eq:interpretation}
\begin{equation}\label{eq:update}
\mathbf U_{K}^{n+1}=\mathbf U_{K}^n-\dt M_K^{-1}\mathbf F_K
\end{equation}
One possible scheme is then: \emph{change nothing for the internal DoFs (cell average) $\xbar \mbu_K$} thanks to \eqref{equation:dg}-\eqref{equation:average}, and for the boundary DoFs, inspired by the residual distribution method, we do:
\begin{equation}\label{eq:projection}
\mathbf U_{\sigma}^{n+1}=\sum_{K, \sigma\in K}{\omega}_{K,\sigma}\mathbf U_{K,\sigma}^{n+1}
\end{equation}
\end{subequations}
with a corresponding weight ${\omega}_{K,\sigma}$ that are assumed to be positive in the scalar case and a matrix in the system case. For consistency reasons, it has to satisfy
\begin{equation*}
  \sum_{K, {\sigma}\in K}\omega_{K,{\sigma}}=1 \text{ (resp. } \Id)
\end{equation*}
This amounts to write, from \eqref{eq:Project} and \eqref{eq:interpretation},
\begin{equation*}
  \mathbf U_\sigma^{n+1}= \mathbf U_\sigma^n-\Delta t \sum_{K, \sigma\in K} \frac{\omega_{K,{\sigma}}}{\vert K\vert}\mathcal{P}\mathbf F_{K,{\sigma}},
\end{equation*}
which can be further written as a Residual Distribution method:
\begin{equation}\label{eq:HO_RD} 
  \mathbf U_\sigma^{n+1}= \mathbf U_\sigma^n-\Delta t \sum_{K, \sigma\in K} \red{\Phi_{K,\sigma}^{\text{HO}}},\quad \red{\Phi_{K,\sigma}^{\text{HO}}}:=\omega_{K,\sigma}\Phi_{K,\sigma}~\mbox{with}~\Phi_{K,\sigma}:=\frac{1}{\vert K\vert}\mathcal{P}\mathbf F_{K,\sigma}.
\end{equation}
The next question is naturally how to choose the weights $\omega_{K,\sigma}$ so that the scheme is stable and consistent? An idea comes to our mind is to find the weights $\omega_{K,\sigma}$ so that we recover the upwinding \pampa introduced in \cite{Abgrall2023a,Abgrall2024a,Abgrall2025a}. We, therefore, take
\begin{equation}\label{eq:upwind_weights}
{\omega}_{K,\sigma}=\bigg (\sum_{K, \sigma\in K} \big (\text{sign}\big(\mbJ(\mbu_\sigma)\cdot\mbn_\sigma^K\big) +\varepsilon\Id_k\big )\bigg)^{-1}\; \bigg ( \text{sign}\big (\mbJ(\mbu_\sigma)\cdot\mbn_\sigma^K\big )+\varepsilon \Id_k\bigg )
\end{equation}
where $\mbJ=(A,B)$ and $A,B\in M_k(\R)$ are the Jacobians of the flux in the $x$- and $y$-direction, respectively. $\mbJ(\mbu_\sigma)\cdot\mbn_\sigma^K=An_x+Bn_y$ with $\mbn_{\sigma}^K=(n_x,n_y)$ being the normal plotted in Figure \ref{fig:1}. For any $\sigma$ over the element $K$, the following hold:
\begin{itemize}
  \item If $\sigma$ is a vertex, $\mbn_\sigma$ is the inward normal vector to the edge opposite to $\sigma$. For example,  $\mbn_1$. 
  \item If $\sigma$ is a midpoint, $\mbn_\sigma$ is the outward normal vector to the edge on which $\sigma$ is sitting. For example, $\mbn_4$. 
\end{itemize}
Finally, $\varepsilon$ is a  positive parameter  to avoid inversion problems. This can occur for example in the case where $\mbf(\mbu)=\mba\mbu$, $\sigma$ is the mid-point of an edge parallel to $\mba$. In that case up-winding makes no sense and the weight ${\omega}_{K,\sigma}$ is equal to $\tfrac{1}{2}$. Another case where it may occur is at the boundaries, for example at supersonic points: if we do not treat carefully the boundary condition, the matrices $\text{sign}\big(\mbJ(\mbu_\sigma)\cdot\mbn_\sigma^K\big)$ may well equal to zero. So considering $\varepsilon\neq 0$ is very useful. To make a parallel with finite differences, the formula \eqref{eq:upwind_weights} is a bit like introducing an upwind bias. In the numerical experiments, we have chosen $\varepsilon=0$ or $\varepsilon=\tfrac{\vert K\vert}{2}$. 

\begin{remark}[Other possible schemes]
Other possibilities of weights can be considered, like taking an arithmetic average. However, coupled with a 3rd order SSP RK time stepping, this scheme is unstable. In the solid advection of a Gaussian spot, this  instability manifests itself after several hundred of time steps. If the velocity field is constant, we do not see this instability, probably because not enough time steps have been done. We have tried to modify the \dg terms in various ways, all depending on parameters that seem tricky to find. For this reason, we have no longer considered this solution, which works beautifully in one dimension, and have the advantage of simplicity.
\end{remark}

\blue{
\begin{remark}[Equivalence with the schemes of \cite{Abgrall2025a}]
One can easily show that in the case of a linear problem, the current version of the scheme and the one developed in \cite{Abgrall2025a} are exactly the same. In the non linear case, this equivalent is lost. In a future paper, we will discuss and evaluate these differences.
\end{remark}
}
\blue{
\begin{remark}
We note that the upwind \pampa scheme has been introduced and successfully tested in \cite{Abgrall2024a,Abgrall2025a} within the framework of finite difference methods and Residual Distribution schemes, in particular the Low Diffusion Advection (LDA) scheme. One of the main motivations for reinterpreting the upwind \pampa scheme within the conventional \dg framework is to enable a more systematic treatment of boundary conditions, without the need to introduce ghost cells. On the other hand, the above derivations show the connection of the standard third-order \pampa scheme with \dg schemes. In a way, the \pampa scheme is a continuous discontinuous Galerkin scheme: it employs a globally continuous approximation commonly used in Continuous Galerkin (CG) method; when using a Runge--Kutta (RK) time procedure, for each RK cycle, it applies one step of \dg followed by a projection onto globally continuous approximations. This new reinterpretation demonstrates the nature of \pampa scheme between CG and \dg.
\end{remark}}

\subsection{Low-order \pampa scheme}
\begin{figure}[ht!]
\centerline{\includegraphics[trim=0.01cm 0.001cm 0.01cm 0.001cm,clip,width=4.0cm]{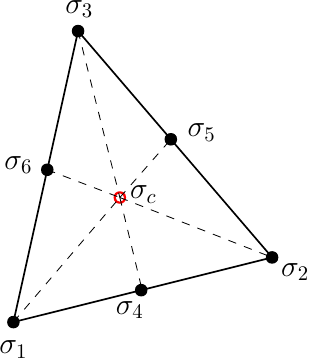}}
\caption{Geometry for the first-order scheme.\label{fig:nscheme} }
\end{figure}

It has been shown that the standard high-order \pampa scheme is neither BP nor OE; see, e.g., \cite{Abgrall2023a,Abgrall2025a}. A remedy is to employ its low-order counterpart. Whenever the numerical solution violates either the BP or OE criterion, the low-order scheme is activated. Here, we briefly recall the key components of the low-order schemes introduced in \cite[Section 3.2]{Abgrall2023a}. For updating the internal DoF $\xbar{\mbu}_K$, or more precisely, to compute the numerical fluxes across each edge of element $K$, we simply consider a standard finite volume discretization. For instance, in our simulations, the local Lax--Friedrichs numerical flux is evaluated using the interior and exterior values of the internal DoFs associated with element $K$ (so we do not use the boundary DoFs of each element). For the evolution of boundary DoFs, we employ a sub-triangulation based variant of the local Lax--Friedrichs method. To define the sub-elements, we first (temporarily) denote the element centroid by $\sigma_c$. Then, we define the six sub-elements by the triplets of vertices as
\begin{equation*}
\begin{aligned}
 &T^K_1=\{\sigma_1,\sigma_c,\sigma_4\},\quad T^K_2=\{\sigma_4,\sigma_c,\sigma_2\},\quad
 T^K_3=\{\sigma_2,\sigma_c,\sigma_5\},\\
 &T^K_4=\{\sigma_5,\sigma_c,\sigma_3\},\quad T^K_5=\{\sigma_3,\sigma_c,\sigma_6\},\quad
 T^K_6=\{\sigma_6,\sigma_c,\sigma_1\}.
 \end{aligned}
\end{equation*}
Finally, within each sub-element, we can obtain inward normals $\mbn_{\sigma_i}^{T_i^K}$ using the three vertices of $T_i^K$. The geometric configuration and the corresponding numbering are illustrated in Figure~\ref{fig:nscheme}.

The update of $\mbu_\sigma$ is done with \eqref{eq:HO_RD} but with the low-order residuals $\Phi_{K,\sigma}^{\text{LO}}$ given by
\begin{equation}\label{eq:LO_RD}
\Phi_{K,\sigma}^{\text{LO}}=\frac{1}{\vert C_\sigma\vert }\sum_{T^K_i, \sigma\in T^K_i}\Phi_{\sigma}^{T^K_i}(\mbu_{\mathrm h}),
\end{equation}
where $\vert C_\sigma\vert=\sum_{K, \sigma\in K}\sum_{T^K_i, \sigma\in T^K_i} \frac{\vert T^K_i\vert}{3}$ is the measure of the dual control volume,
\begin{equation}\label{eq:LO_RD2}
\begin{aligned}
  \Phi_\sigma^{T^K_i}&=\frac{1}{3}\int_{T^K_i} \mbJ(\mbu_{\mathrm h})\cdot \nabla \mbu_{\mathrm h} \; {\rm d}\mbx+\alpha_{T^K_i}\big (\mbu_\sigma-\overline{\mbu}_{T^K_i}\big )\\
  &\approx\frac{1}{3}\vert T^K_i\vert \mbJ(\xbar{\mbu}_{T^K_i})\nabla \mbu_{\mathrm h}(\red{\overline{\bbx}_{T^K_i}})+\alpha_{T^K_i}\big (\mbu_\sigma-\overline{\mbu}_{T^K_i}\big )
\end{aligned}
\end{equation}
with \red{$\overline{\bbx}_{T^K_i}$ being the arithmetic average of the three vertices of $T^K_i$}, $\overline{\bbu}_{T^K_i}$ being the arithmetic average of the $\bbu$'s at the three vertices of $T^K_i$ (hence we use the average value here), and 
\begin{equation}\label{eq:LO_RD3}
  \alpha_{T^K_i}\geq \max\limits_{\substack{\sigma\in T^K_i,\\\mbn_{\sigma_i}^{T_i^K} \text{normals of~} T^K_i}}\rho\big ( \mbJ(\mbu_\sigma)\cdot \mbn_{\sigma_i}^{T_i^K}\big ),
\end{equation}
where $\rho(B)$ is the spectral radius of the matrix $B$.

\subsection{Bound Preservation and Oscillation Elimination}
Here, we address how to naturally combine the standard high-order and low-order \pampa schemes. Our goal is to enhance the robustness of the high-order \pampa scheme by incorporating low-order discretizations while preserving its high-order accuracy as much as possible. An effective strategy is to perform a local convex blending of the high-order and low-order fluxes or residuals, as recently developed in \cite{Abgrall2026,Abgrall2024a}.

Let $\hat{\mbf}^{\text{HO}}{\mbn_e}$ denote the high-order numerical flux across an edge $e$ of element $K$, computed using Gaussian quadrature, and let $\hat{\mbf}^{\text{LO}}{\mbn_e}$ represent the corresponding low-order flux. The blended flux is then defined as
\begin{equation}\label{eq:eta_e}
  \hat{\mbf}_{\mbn_e}^{\text{blend}}=\hat{\mbf}^{\text{LO}}_{\mbn_e}+\eta_e\big(\hat{\mbf}^{\text{HO}}_{\mbn_e}-\hat{\mbf}^{\text{LO}}_{\mbn_e}\big)=\hat{\mbf}^{\text{LO}}_{\mbn_e}+\eta_e\Delta\hat{\mbf}_{\mbn_e}.
\end{equation}
Similarly, the blended residual for a boundary DoF $\sigma$ is given by
\begin{equation}\label{eq:eta_sigma}
  \red{\Phi_{K,\sigma}^{\text{blend}}}=\Phi_{K,\sigma}^{\text{LO}}+\eta_\sigma\big(\Phi_{K,\sigma}^{\text{HO}}-\Phi_{K,\sigma}^{\text{LO}}\big)=\Phi_{K,\sigma}^{\text{LO}}+\eta_\sigma\Delta\Phi_{K,\sigma},
\end{equation}
where the high-order residual $\Phi_{K,\sigma}^{\text{HO}}$ is given by \eqref{eq:HO_RD}--\eqref{eq:upwind_weights} and the low-order ones in \eqref{eq:LO_RD}--\eqref{eq:LO_RD3}. By substituting these blended fluxes and residuals into the evolution equations for the internal and boundary DoFs, respectively, the updated solution can be expressed as a convex combination of the solution from the previous time step and several intermediate states. We, therefore, require that all intermediate states remain within the invariant domain $\DD$, from which sufficient conditions on the blending parameters $\eta_e$ and $\eta_\sigma$ can be derived. For brevity, we omit the derivations of the explicit expressions of these BP parameters and refer the readers to \cite{Abgrall2024a}. 
It is worth emphasizing that when this technique is applied to the Euler equations, we have derived an explicit and robust computation of the BP parameters using the GQL representation \eqref{eq:2}. This effective and efficient approach for determining the BP parameters is not restricted to the \pampa scheme, it also shows the potential for other numerical methods; see, e.g., \cite{Wissocq2025,Kuzmin2025}.    

\paragraph*{Convex Oscillation Eliminating} It is important to note that being BP does not mean that we are non-oscillatory. To mitigate nonphysical waves that may arise near strong discontinuities, we revisit the blending weights of the high- and low-order fluxes and residuals by incorporating information from the jumps of high-order derivatives at element edges. Inspired by the strategy in \cite{Abgrall2025c} and following ideas from \cite{Peng2025,Lu2021,Abgrall2025b}, we introduce an element-based convex OE parameter $\theta_K \in [0,1]$ defined as
\begin{equation}\label{eq:OE_param}
    \theta_K=\exp\Big(-\frac{1}{N_K}\sum_{e\subset\partial K}\frac{\alpha_e\dt_n}{\ell_{e,K}}\sigma_{e,K}(\mbu_{\mathrm h})\Big),
\end{equation}
where $N_K$ is the cardinality of the neighborhood element set of triangle $K$, namely,
\begin{equation*}
  N_K=\#\underline{V}_K=\#\{K':\partial K\cap \partial K'=e\subset\partial K\};  
\end{equation*}
for triangular meshes, $N_K=3$. Here, $\dt_n$ is the $n$-th adaptive time step computed using a suitable CFL condition, and $\ell_{e,K}$ is defined as
\begin{equation*}
    \ell_{e,K}=\sup_{\mbx\in K}\rm{dist}~(\mbx,e).
\end{equation*}
The damping parameter $\sigma_{e,K}(\mbu_{\mathrm h})$ plays a critical role in balancing accuracy and stability (i.e., effectively suppress spurious oscillations here). It is defined as
\begin{equation}
    \sigma_{e,K}(\mbu_{\mathrm h})= \sigma_{e,K}\big(\text{RI}(\mbu_{\mathrm h})\big)=\max_{1\leq i\leq p}\sigma_{e,K}\big(\widehat\mbu_{\mathrm h}^{(i)}\big), 
\end{equation}
where $\widehat\mbu_{\mathrm h}=\text{RI}(\mbu_{\mathrm h})$ denotes a rotational invariant reformulation:
\begin{itemize}
    \item In the scalar case, $\text{RI}(\mbu_{\mathrm h})=\mbu_{\mathrm h}$: we do nothing
    \item In the Euler case, 
according to \cite{Ding2025}, we introduce the normal and tangential components of the momentum:
\begin{equation}
\begin{split}
  \widehat{\mbu}^{(1)}_{\mathrm h}&=\mbu^{(1)}_{\mathrm h},~\widehat{\mbu}^{(4)}_{\mathrm h}=\mbu^{(4)}_{\mathrm h}\\
   \widehat{\mbu}^{(2)}_{\mathrm h}&=\mbu^{(2)}_{\mathrm h}\cdot\mbn_e^{(1)}+\mbu^{(3)}_{\mathrm h}\cdot\mbn_e^{(2)},~
  \widehat{\mbu}^{(3)}_{\mathrm h}=-\mbu^{(2)}_{\mathrm h}\cdot\mbn_e^{(2)}+\mbu^{(3)}_{\mathrm h}\cdot\mbn_e^{(1)}
  .
  \end{split}
\end{equation}
\end{itemize}
Then, the component-wise quantity is given by
\begin{equation}\label{eq:damping}
\sigma_{e,K}(\widehat\mbu_{\mathrm h}^{(i)})=\left\{\begin{aligned}
&0, &&\mbox{if}~\widehat\mbu_{\mathrm h}^{(i)}\equiv\langle\widehat\mbu_{\mathrm h}^{(i)}\rangle_\Omega\\
&\sum_{\kappa=1}^{k=2}C_{\kappa}\ell_{e,K}^\kappa\sum_{\vert\bm{\alpha}\vert=\kappa}\dfrac{\dfrac{1}{\vert e\vert}\int_e\vert \llbracket \partial^{\bm{\alpha}}\widehat\mbu^{(i)}_{\mathrm h}\rrbracket_e\vert\;{\mathrm d}\ell}{\Vert \mbu_{\mathrm h}^{(i)}-\langle \mbu_{\mathrm h}^{(i)}\rangle_\Omega\Vert_{L^\infty(\Omega)}},&&\mbox{otherwise.}
\end{aligned}\right.
\end{equation}
Here, $C_{\kappa}$ is a problem-independent parameter and we simply take $C_{\kappa}=1$ in this paper, $\vert e\vert$ denotes the length of the edge $e\subset\partial K$, $\Vert\mbu_{\mathrm h}^{(i)}-\langle\mbu_{\mathrm h}^{(i)}\rangle_\Omega\Vert_{L^\infty(\Omega)}$ is evaluated at all DoFs and the integration on $e$ is approximated by the Gauss quadrature rule and $\langle\mbu_{\mathrm h}^{(i)}\rangle_\Omega=\frac{1}{\vert\Omega\vert}\int_\Omega \mbu_{\mathrm h}^{(i)}\;\mathrm d\mbx$. $\bm\alpha=(\alpha_1,\ldots,\alpha_d)$ is the multi-index with $\vert\bm\alpha\vert=\sum_{j=1}^d\alpha_j$ and the operator $\partial^{\bm\alpha}$ is defined as 
\begin{equation*}
    \partial^{\bm\alpha}=\dfrac{\partial^{\vert\bm\alpha\vert}}{\partial x_1^{\alpha_1}\cdots\partial x_d^{\alpha_d}}.
\end{equation*}
In the numerical implementation, to mitigate the effects of round-off errors, the condition $\widehat\mbu_{\mathrm h}^{(i)}\equiv\langle\widehat\mbu_{\mathrm h}^{(i)}\rangle_\Omega$ is replaced with
\begin{equation*}
  \Vert\mbu_{\mathrm h}^{(i)}-\langle\mbu_{\mathrm h}^{(i)}\rangle_\Omega\Vert_{L^\infty}\leq10^{-12}\max\big(1,\langle\mbu_{\mathrm h}^{(i)}\rangle_\Omega\big).
\end{equation*}
In the end, the convex BP OE parameters for updating internal (resp. boundary) DoFs are defined as the minimal values between the BP parameters $\eta_e$ (resp. $\eta_\sigma$) and the OE parameters given in \eqref{eq:OE_param}. \red{Namely, we take
\begin{equation*}
  \eta_e^{\text{BP OE}}=\min\big(\eta_e^{\text{BP}},\theta_K,\theta_{K'}\big), \quad \partial K\cap \partial K'=e
\end{equation*}
and
\begin{equation*}
  \eta_\sigma^{\text{BP OE}}=\min\big(\eta_\sigma^{\text{BP}}, \theta_K\big),\quad \sigma\in K,
\end{equation*}
where $\eta_e^{\text{BP}}$ and $\eta_\sigma^{\text{BP}}$ are given in \cite{Abgrall2024a}.}

\section{Boundary Conditions}
Handling boundary conditions in \pampa scheme is somewhat subtle, particularly on unstructured meshes, because part of the algorithm is inspired by finite-difference ideas. Although this can be done, as shown in \cite{Abgrall2024a}, we now adopt a more robust and systematic approach based on the \dg formulation. 

Consider the problem \eqref{eq:1} in a domain $\Omega$ with the boundary condition
\begin{equation}\label{eq:hyper_bc}
\big ( \mbJ\cdot \mbn\big )^-\big (\mbu(\mbx,t)-\mbu_b(\mbx, t)\big )=0, \quad t>0
\end{equation}
where $\mbn$ denotes the unit outward normal at $\mbx\in\partial\Omega$, $\mbu_b\in\mathbb{R}^m$ is a prescribed boundary state, and $\mbJ$ is the Jacobian of $\mbf$ with respect to $\mbu$. Recall that the problem is hyperbolic. We assume that $\partial\Omega$ is smooth; see \cite{Ern2006} for more details. The most natural modification is to replace \eqref{eq:dg_load} by
\begin{equation}\label{eq:dg boundary}
\begin{aligned}
F_\sigma&=
    -\int_K \nabla \varphi_\sigma \cdot \mbf(\mbu_{\mathrm h})\; \mathrm d\mbx
    +\int_{\partial K}\varphi_\sigma\hat{\mbf}_\mbn(\mbu_{\mathrm h},\mbu_{\mathrm h}^-) \; \mathrm d\ell,\\
F_\mu&=-\int_{K} \nabla \varphi_{7} \cdot \mbf(\mbu_{\mathrm h})\; \mathrm d\mbx+\int_{\partial K}\varphi_{7}\hat{\mbf}_\mbn(\mbu_{\mathrm h},\mbu_{\mathrm h}^-) \; \mathrm d\ell.
\end{aligned}
\end{equation}
where $\hat{\mbf}_\mbn$ is a numerical flux consistent with \eqref{eq:hyper_bc} on the edge $e_b=\partial K\cap \partial \Omega$ for which we set $\mbu_{\mathrm h}^-=\mbu_b$. For the cell average, this implies the addition of
\begin{equation}\label{eq:BCmodify_ave}
  \int_{e_b}\big ( \hat{\mbf}_\mbn(\mbu_{\mathrm h}, \mbu_{\mathrm h}^-)-\mbf(\mbu_{\mathrm h})\cdot \mbn\big )\; \mathrm d\ell,
\end{equation}
and for the point values $\sigma\in K$ that shares an edge on $\partial \Omega$, we get
\begin{equation*}
 \Phi_{\sigma,K}=-\int_K\mbf(\mbu_{\mathrm h})\nabla\theta_\sigma\;\mathrm d\mbx+\int_{\partial K}\theta_\sigma\hat{\mbf}_{\mbn}(\mbu_{\mathrm h},\mbu_{\mathrm h}^-)\; \mathrm d\ell,
\end{equation*}
i.e., if $e$ is the edge in $\partial\Omega$,
\begin{equation*}
  \Phi_{\sigma,K}=-\int_K\mbf(\mbu_{\mathrm h})\nabla\theta_\sigma\;\mathrm d\mbx+\int_{\partial K}\theta_\sigma\mbf(\mbu_{\mathrm h})\cdot {\mbn}\; \mathrm d\ell+\int_e \theta_\sigma\big ( \hat{\mbf}_\mbn(\mbu_{\mathrm h}, \mbu_{\mathrm h}^-)-\mbf(\mbu_{\mathrm h})\cdot {\mbn}\big ) \; \mathrm d\ell,
\end{equation*}
This shows that we need to add 
\begin{equation}\label{eq:BCmodify_point}
  \int_e \theta_\sigma\big ( \hat{\mbf}_\mbn(\mbu_{\mathrm h}, \mbu_{\mathrm h}^-)-\mbf(\mbu_{\mathrm h})\cdot {\mbn}\big ) \; \mathrm d\ell
\end{equation}
to the already computed residuals. From this expression, we see from \eqref{eq:dg_load} that \textit{all} the DoFs of $K$ must be modified.

In this paper, we have considered two types of boundary conditions:
\begin{itemize}
    \item Inflow/outflow conditions. This occurs in the scalar and system case. The state $\mbu_h^-$ is the condition at infinity. We have taken the Steger-Warming flux for its simplicity. Another choice could be the Local Lax-Friedrich flux.
    \item Wall condition. This occurs only for the Euler equations. Here the state $\mbu_h^-$ is obtained from $\mbu_h$ by symmetry with respect to the normal: same density, total energy, and the momentum is symmetrized. Then we take the local Lax-Friedrich flux.
\end{itemize}
\begin{remark}
We first remark that in the computation of the first parts of the integrals in \eqref{eq:BCmodify_ave} and \eqref{eq:BCmodify_point} (those applied the numerical fluxes $\hat{\mbf}{\mbn}$), we also adopt Gaussian quadrature rules. At each quadrature point, the values of $\mbu_{\mathrm h}$ and $\mbu_{\mathrm h}^-$ must remain in the invariant domain $\DD$, which is ensured by the local scaling BP limiter introduced in \cite{Zhang2010,Zhang2010a}. We also notice that these quadrature-based computations may violate the BP property. A simple remedy is to evaluate the above integrals using only the internal DoF $\xbar{\mbu}_K$ and its prescribed boundary state $\xbar{\mbu}_K^-$ in the numerical fluxes. A more accurate alternative is to blend the Gaussian quadrature version with this simplified variant. In the numerical experiments presented in the next section, we employ the Gaussian quadrature version and do not observe any violations of the BP property.
\end{remark}
\section{Numerical Examples}
In this section, we present a comprehensive suite of benchmark and challenging numerical experiments to demonstrate the high-order accuracy, robustness and effectiveness of the proposed BP OE \pampa scheme on unstructured triangular meshes. In all examples, we employ the \red{three-stage, third-order} strong-stability-preserving Runge--Kutta (SSPRK 33) time discretization method \red{and at each RK stage, we therefore need to apply the projection step and the BP OE limiter}. The unstructured triangular meshes are generated using GMSH \cite{Geuzaine2009}. Unless otherwise specified, the adiabatic index for Euler equations is set to $\gamma=1.4$ and the CFL number $0.2$.

\phantomsection
\subsection{Scalar case}
We begin with the scalar case, 
\begin{equation*}
\dpar{u}{t}+\mba(\mbx)\cdot\nabla u=0,
\end{equation*}
where $\mba(\mbx)=\big (a^{(1)},a^{(2)}\big )$ represents the velocity field that the quantity $u$ is moving with. Three test cases are taken into account here. The boundary conditions for all the scalar problems are chosen such that the solution does not change in a neighborhood of the boundary.

\subsubsection*{Example 1---Transport Problem}
In the first example, we consider $\mba=(-1,-1)$ and the initial
condition 
\begin{equation*}
    u_0(\mbx)=\mathrm e^{-\frac{1}{4}\Vert \mbx-\mbx_0\Vert^2}, ~~\mbx_0=(15, 15).
\end{equation*}
over a computational domain $[-20,20]^2$.
The solution is checked at a long final time $T=10$ with the exact one. The discrete $L^1$-, $L^2$-, and $L^{\infty}$-errors of the internal and boundary DoFs are shown in
Tables \ref{tab:1}. As one can clearly see that, the expected third-order accuracy is achieved.
\begin{table}[ht!]
 \caption{Example 1: Transport problem. Error analysis and experimental convergence rate.\label{tab:1}}
\begin{center}
\begin{tabular}{||c||cc||cc||}
\hline
&\multicolumn{2}{c||}{Internal DoF}&\multicolumn{2}{c||}{Boundary DoFs}\\
\hline
$h$& $L^1$-error &rate &$L^1$-error &rate\\\hline
 $1.1141$& $0.1348\,10^{-2}$& -&$0.9756\,10^{-3}$&-\\
 $0.5570$ &$0.3068\,10^{-3}$& $2.135$&$0.2087\,10^{-3}$ &$2.225$\\
 $0.2785$ &$0.4473\,10^{-4}$ &$2.778$&$0.3005\,10^{-4}$& $2.796$\\
 $0.1542$ &$0.8360\,10^{-5}$ &$2.837$ &$0.5614\,10^{-5}$& $2.838$\\
 \hline
 \multicolumn{5}{||c||}{}\\
 \hline
$h$& $L^2$-error &rate &$L^2$-error &rate\\  
\hline
$1.114$& $0.7618\,10^{-2}$& -&$0.6669\,10^{-2}$&-\\
$0.5570$ &$0.1895\,10^{-2} $&$2.007$&$0.1585\,10^{-2}$& $2.073$\\
$ 0.2785$& $0.2871\,10^{-3}$ &$2.723$&$0.2368\,10^{-3}$&$ 2.743$\\
 $0.1542$ &$0.5367\,10^{-4}$ &$2.836$&$0.4413\,10^{-4}$& $2.842$\\
 \hline
 \multicolumn{5}{||c||}{}\\
 \hline
 $h$& $L^\infty$-error &rate &$L^\infty$-error &rate\\  
 \hline
 $1.114$& $0.1251$& -&$0.1335$& -\\
 $0.5570$& $0.3424\,10^{-1}$& $1.870$& $0.3597\,10^{-1}$& $1.892$\\
 $0.2785$& $0.5518\,10^{-2}$ &$2.634$&$0.5608\,10^{-2}$ &$2.681$\\
 $0.1542$ &$0.1045\,10^{-2}$ &$2.815$&$0.1064\,10^{-2 }$&$2.811$\\
 \hline
 \end{tabular}
 \end{center}
 \end{table}

\subsubsection*{Example 2---Zalesak Problem}
In the second example, taken from \cite{Abgrall2025a}, we consider  Zalesak's problem, which involves the solid body rotation of a notched disk. Here the domain is $[0,1]^2$ and the rotation field is $\mba(\mbx)$ also defined to be a rotation with respect to $(x_0,y_0)=(0.5,0.5)$ and angular speed of $2\pi$.
The initial condition is 
\begin{equation*}
 u(\mbx,0)=\left \{
   \begin{array}{ll}
   0.25\big ( 1+\cos(\frac{\pi r_{1}}{0.15})) & \text{ if }r_{1}=\Vert\mbx-(0.25, 0.5)\Vert\leq 0.15,\\
  1-\frac{r_2}{0.15}&\text{ if } r_{2}=\Vert\mbx-(0.5,0.25)\Vert\leq 0.15, \\
  1&\text{ if } r_{3}=\Vert\mbx-(0.5,0.75)\Vert\leq 0.15,
  \\
  0& \text{ if }|x-0.5| \leq 0.025 ~\text{ and }~y\in[0.6,0.85]. 
  \end{array}\right.   
\end{equation*}

\begin{figure}[ht!]
\centerline{\subfigure  [Internal DoF $\xbar u_K$, coarse mesh]{\includegraphics[trim=1.8cm 2.5cm 2.5cm 1.0cm,clip,width=5.5cm]{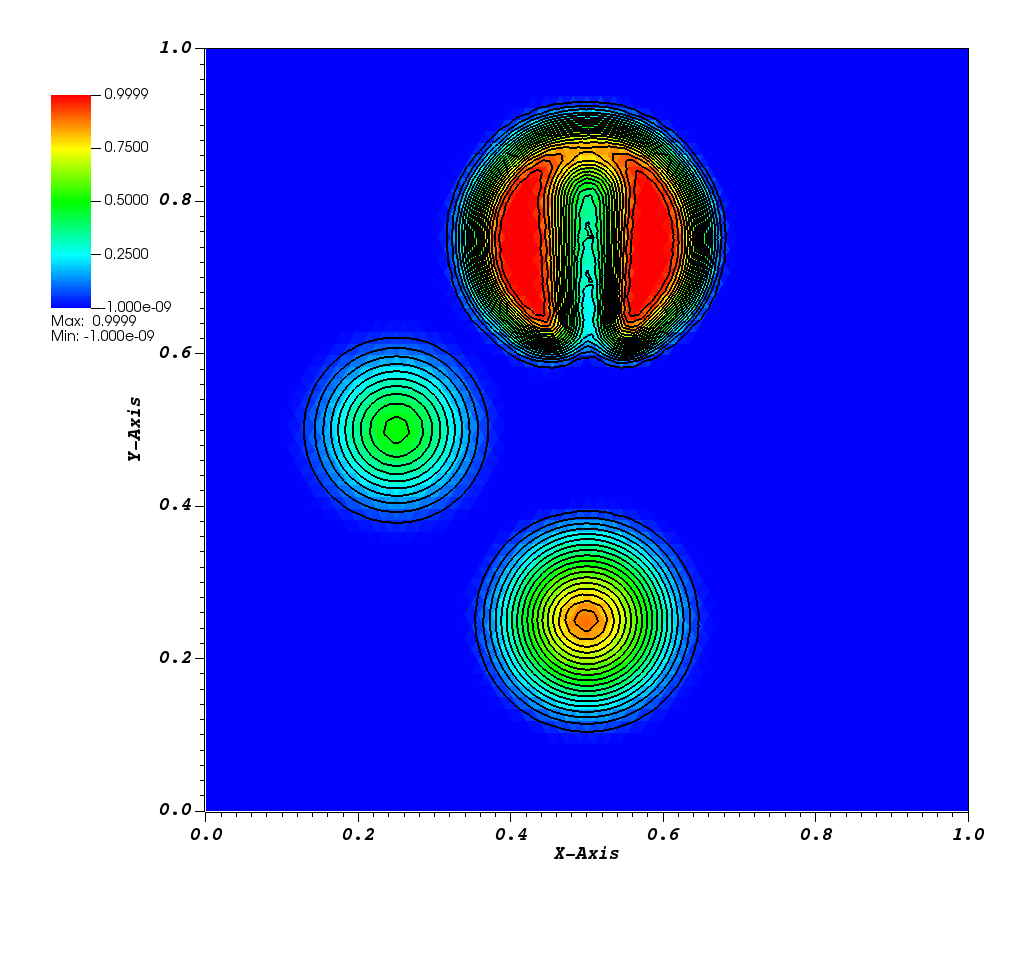}}\hspace*{0.05cm}
	\subfigure[Boundary DoFs $u_\sigma$, coarse mesh]{\includegraphics[trim=1.8cm 2.5cm 2.5cm 1.0cm,clip,width=5.5cm]{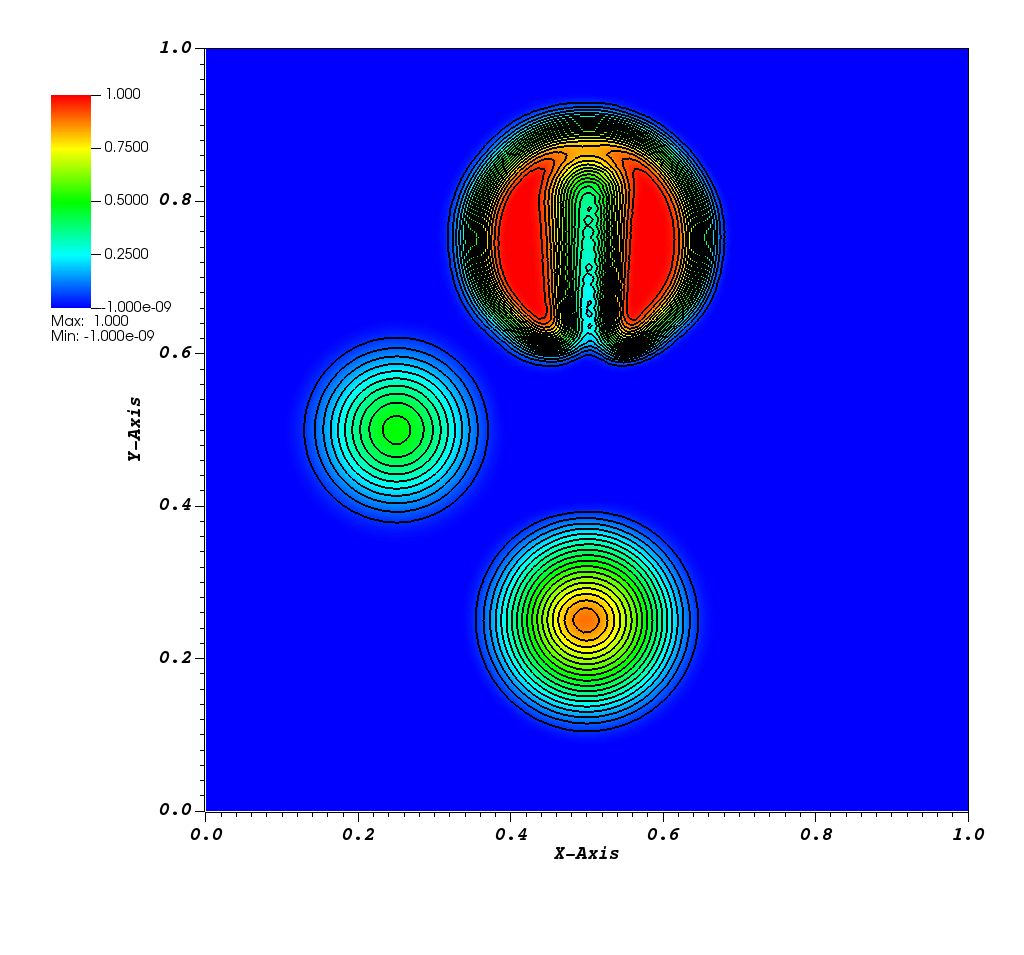}}}
\vskip1pt
    \centerline{\subfigure [Internal DoF $\xbar u_K$, fine mesh]{\includegraphics[trim=1.8cm 2.5cm 2.5cm 1.0cm,clip,width=5.5cm]{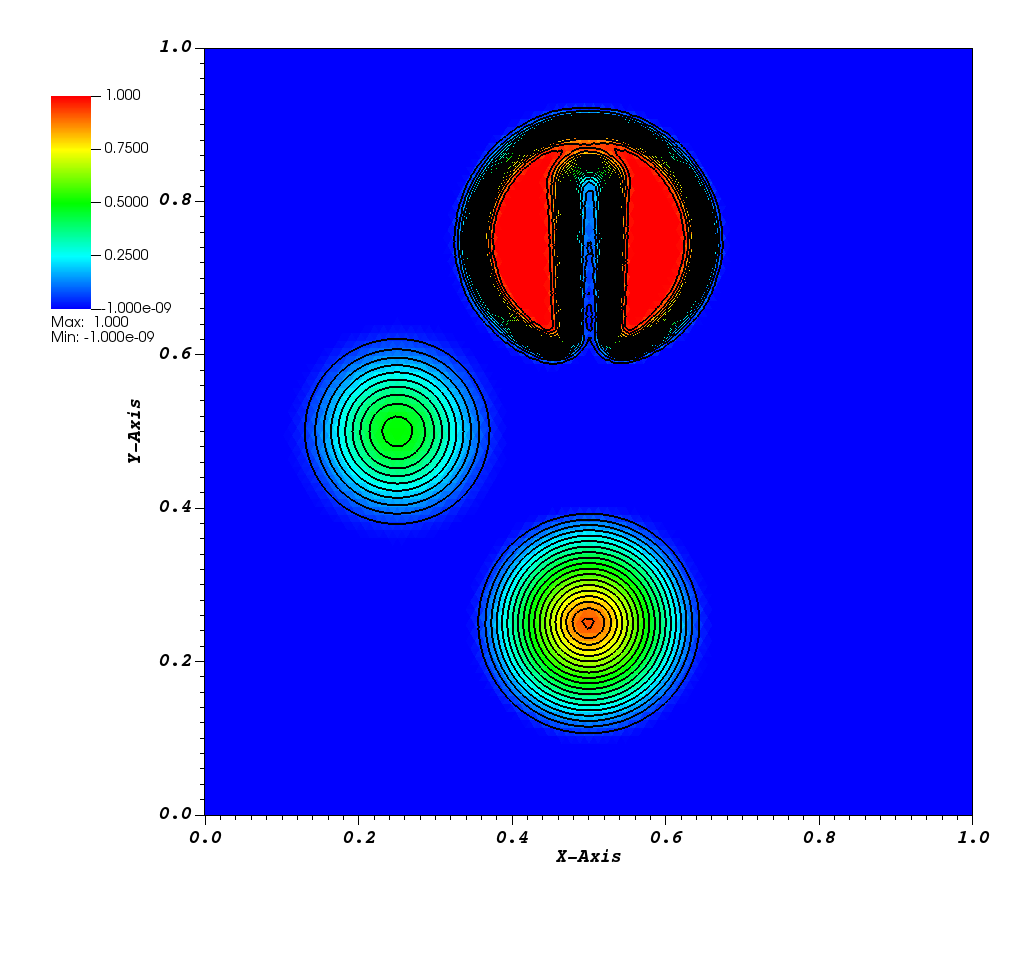}}\hspace*{0.05cm}
	\subfigure[Boundary DoFs $u_\sigma$, fine mesh]{\includegraphics[trim=1.8cm 2.5cm 2.5cm 1.0cm,clip,width=5.5cm]{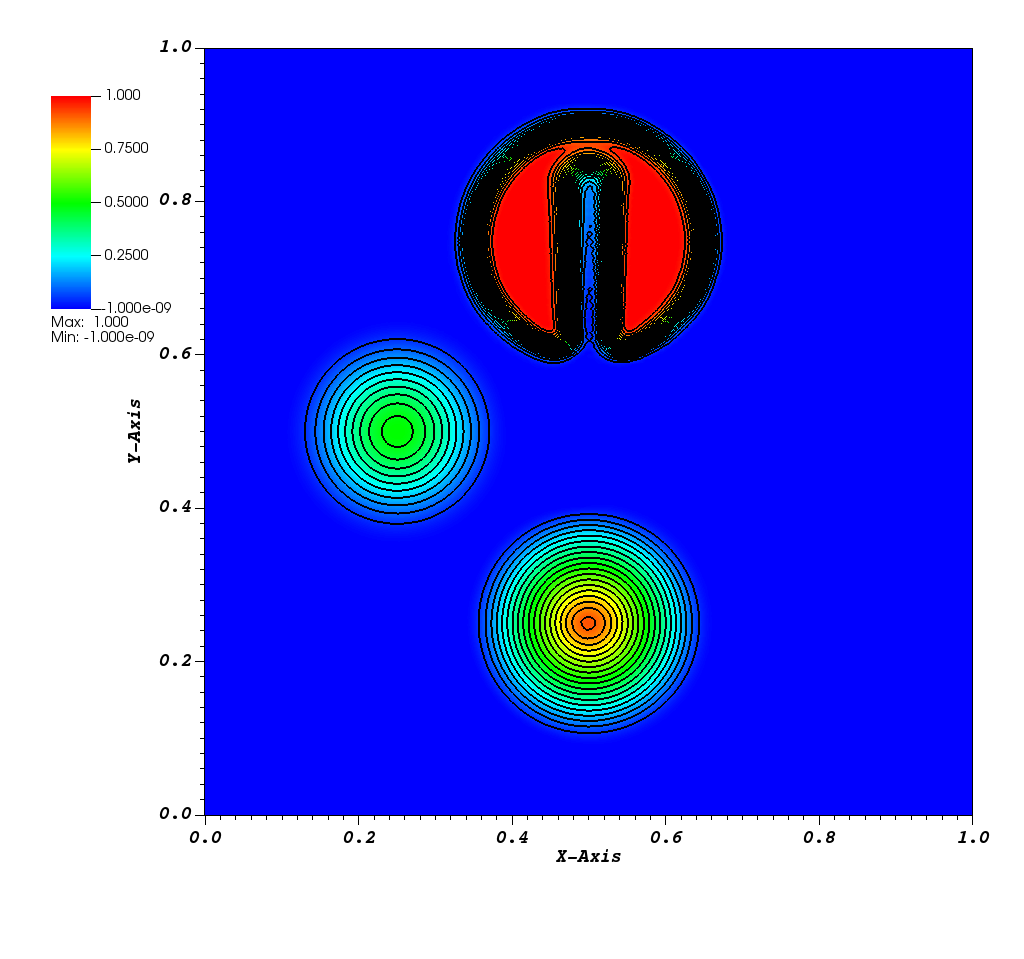}}}
\caption{\sf Example 2: Zalesak problem. Numerical solutions at $t=1$. \blue{Top row: coarse mesh with 7,320 elements; Bottom row: fine mesh with 16,374 elements.} \label{fig: Zalesak}}
\end{figure}

We set the invariant domain $\DD=[-10^{-9},1.0+10^{-9}]$ and compute the numerical solution with a \blue{coarse} mesh that has 7,320 elements and 14,865 points \blue{and a finer mesh with 16,374 elements and 33,085 points}. The results computed by the BP OE \pampa scheme after 1 rotation are displayed in Figure \ref{fig: Zalesak}, where one can see that both the average values and point values can be correctly captured.

\subsubsection*{Example 3---KPP Problem}
In the third example, we consider the so-defined classical KPP (Kurganov--Petrova--Popov) problem, which admits infinitely many weak solutions \cite{Kurganov2007}:
\begin{equation*}
    \dpar{u}{t}+\dpar{\sin(u)}{x}+\dpar{\cos(u)}{y}=0,
\end{equation*}
prescribed in a domain $[-2,2]^2$ with the initial condition
\begin{equation*}
    u(\mbx,0)=\left\{
    \begin{aligned}
    &\frac{7}{2}\pi &&\mbox{if}~\Vert\mbx-(0,0.5)\Vert\leq 1,\\
    &\frac{\pi}{4} &&\mbox{else}.
    \end{aligned}\right.
\end{equation*}
The unique entropy solution at the final time $t=1$ exhibits a two-dimensional rotational wave structure.  A key difficulty of this test is ensuring that the numerical scheme introduces sufficient dissipation so that the approximations converge to the correct entropy solution, and is not too dissipative so that one gets a sharp representation of the discontinuities.

\begin{figure}[ht!]
\centerline{\subfigure[Internal DoF $\xbar u_K$]{\includegraphics[trim=1.8cm 2.5cm 2.5cm 1.0cm,clip,width=5.5cm]{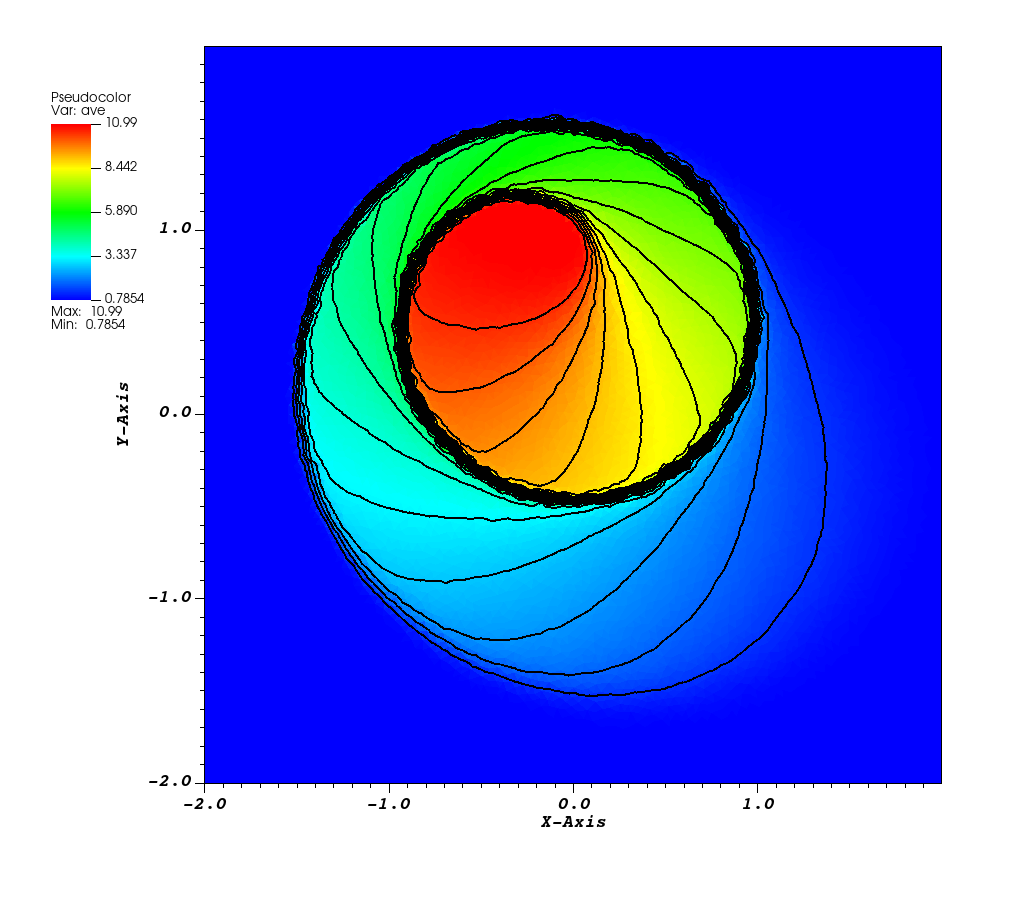}}\hspace*{0.05cm}
	\subfigure[Boundary DoFs $u_\sigma$]{\includegraphics[trim=1.8cm 2.5cm 2.5cm 1.0cm,clip,width=5.5cm]{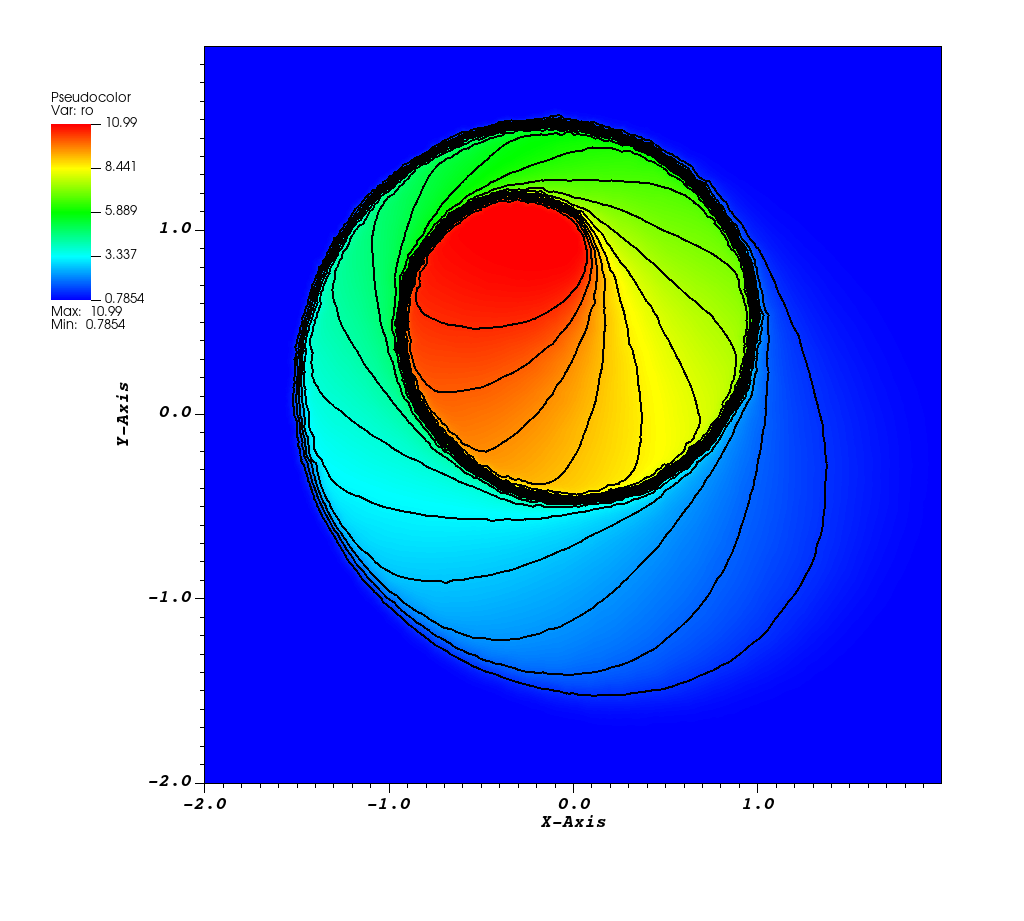}}}
\caption{\sf Example 3: KPP problem. Numerical solutions at $t=1$. \label{fig: KPP}}
\end{figure}

We set the invariant domain $\DD=[-1.0, 100.]$, knowing that the solution must stay in the bound $[\tfrac{\pi}{4},\tfrac{7}{2}\pi]$. It is a way to test if the OE procedure is robust enough not to create undershoot and overshoot in a not controllable manner. We run the computation on an unstructured triangular mesh with 17,200 elements and 34,081 points. Figure \ref{fig: KPP} shows the solutions (cell averages and point values) produced by the BP OE \pampa scheme at the final time. The results are comparable with the published ones; see, e.g., \cite{Abgrall2024a,Abgrall2025a,Kurganov2007}. The two shocks in the northwest direction at around $22^\circ$ are not attached, as expected.
 
\subsection{Euler Equations of Gas Dynamics}
A wide variety of standard and challenging numerical tests for compressible Euler equations are available to evaluate the robustness and accuracy of a proposed numerical scheme. In the following, we select several representative and demanding examples, typical in the literature and sufficiently complex from our perspective, to demonstrate the performance of the proposed BP OE \pampa scheme. In all cases, the invariant domain is such that
$\rho, p\in [10^{-10},10^{10}]$. We simply want that the BP procedure avoids blow-up. In all cases the CFL number is set to $0.3$, and $\gamma=1.4$ except for case 9 where it is set to $5/3$.

\subsubsection*{Example 4---Kurganov--Tadmor Problem}
\begin{figure}[ht!]
\centerline{\subfigure[Internal DoF $\xbar\rho_K$, coarse mesh]{\includegraphics[trim=1.8cm 2.5cm 1.05cm 1.0cm,clip,width=5.5cm]{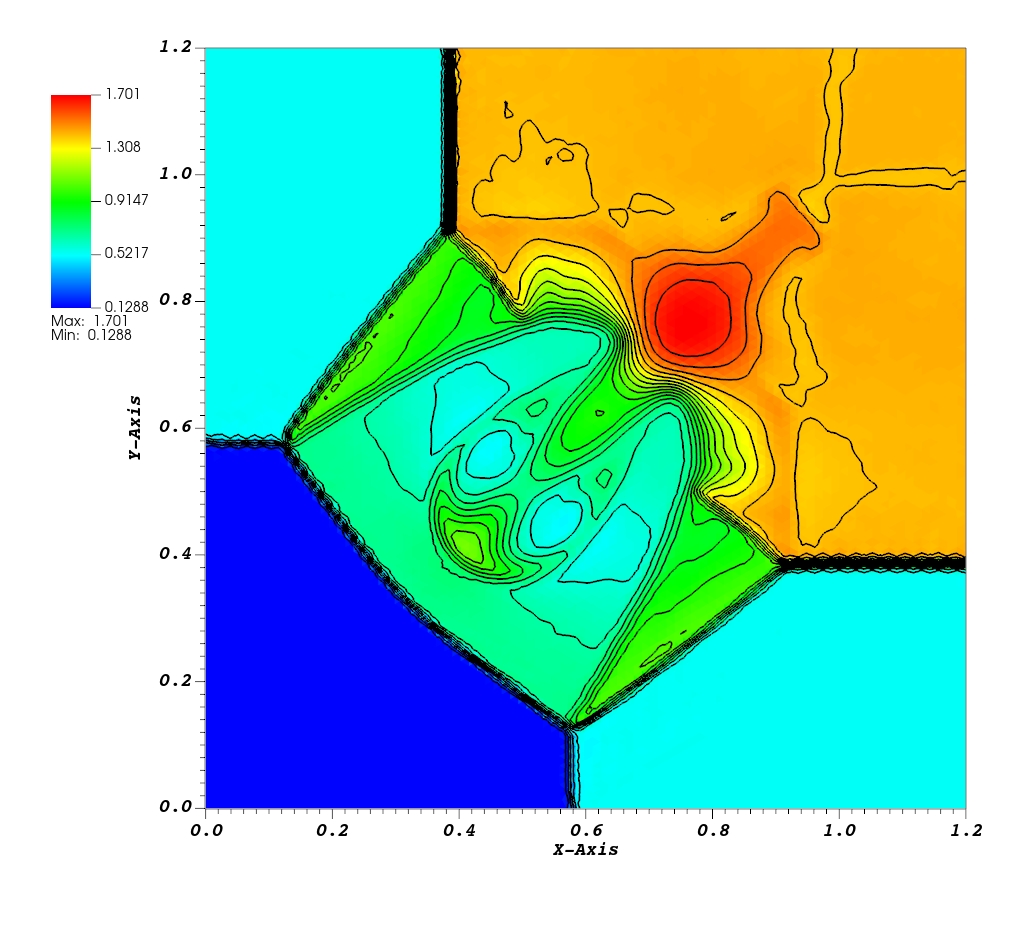}}\hspace*{0.05cm}
	\subfigure[Boundary DoFs $\rho_\sigma$, coarse mesh]{\includegraphics[trim=1.8cm 2.5cm 1.05cm 1.0cm,clip,width=5.5cm]{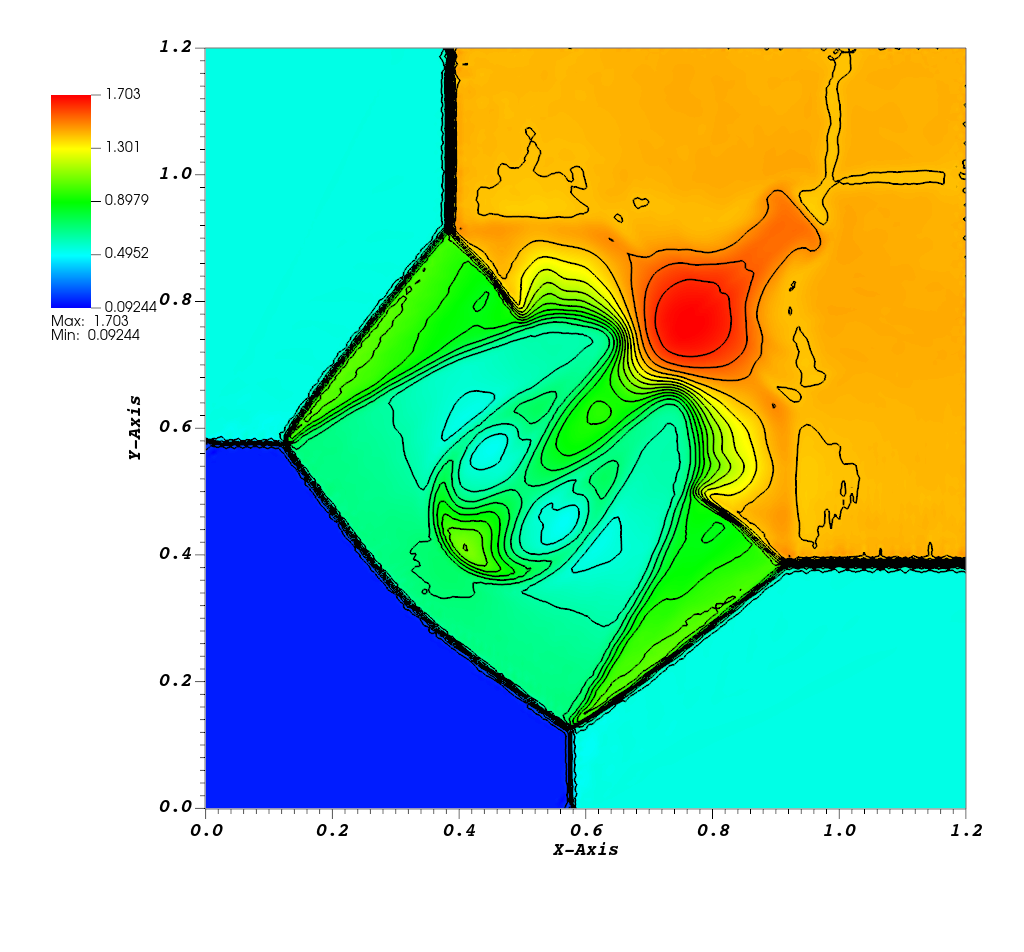}}}
    \vskip1pt
    \centerline{\subfigure[Internal DoF $\xbar\rho_K$, fine mesh]{\includegraphics[trim=1.8cm 2.5cm 1.05cm 1.0cm,clip,width=5.5cm]{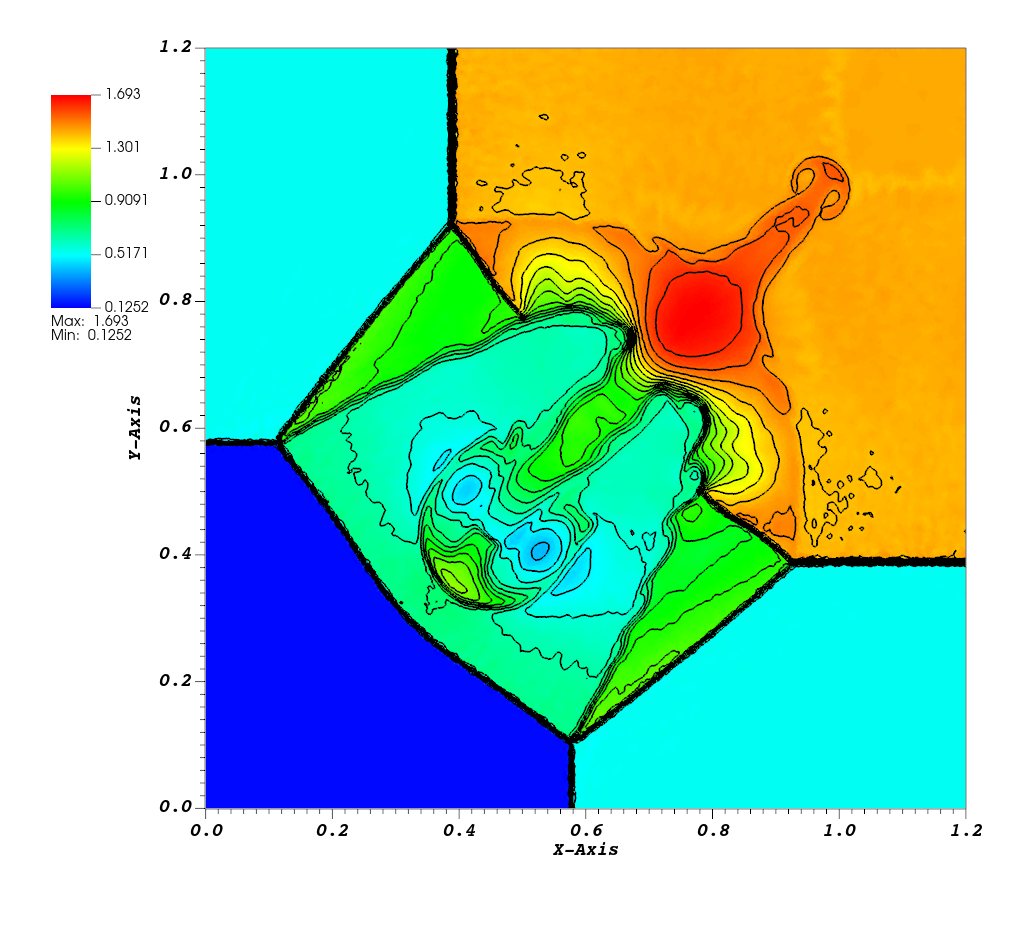}}\hspace*{0.05cm}
	\subfigure[Boundary DoFs $\rho_\sigma$, fine mesh]{\includegraphics[trim=1.8cm 2.5cm 1.05cm 1.0cm,clip,width=5.5cm]{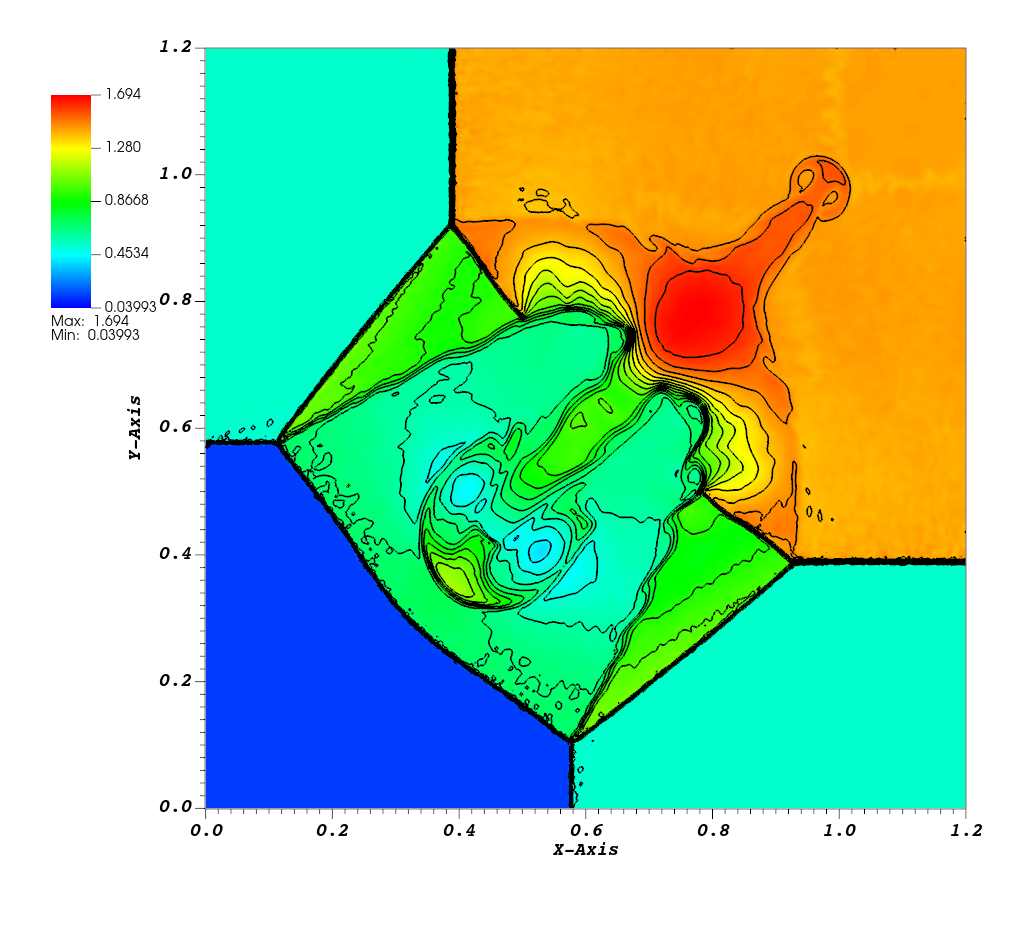}}}
\caption{\sf Example 4: KT problem. Numerical solutions (density) at $t=1$. \blue{Top row: coarse mesh with 15,106 elements; Bottom row: fine mesh with 60,230 elements.}\label{fig:KT}}
\end{figure}

In the fourth example,  which is also known as Configuration 3 of Lax \& Liu in \cite{Lax1998}, we consider the initial condition, 
\begin{equation*}
   (\rho, u,v,p)=\left \{\begin{array}{ll}
(\rho_1,u_1,v_1,p_1)=(1.5, 0, 0, 1.5)& \text{ if } x\geq1\text{ and } y\geq 1,\\
(\rho_2,u_2,v_2,p_2)=(0.5323, 1.206, 0, 0.3) & \text{ if } x\leq 1 \text{ and } y\geq 1,\\
(\rho_3,u_3,v_3,p_3)=(0.138, 1.206, 1.206, 0.029)&\text{ if } x\leq 1\text{ and }y\leq 1,\\
(\rho_4,u_4,v_4,p_4)=(0.5323, 0, 1.206, 0.3) &\text{ if } \red{x\geq1}\text{ and } y\leq 1.
\end{array}\right .
\end{equation*}
prescribed in the domain $[0,1.2]^2$. We compute the numerical solutions on a unstructured triangular mesh with 15,106 elements and 29,893 points \blue{as well as on a finer mesh with 60,230 elements and 119,861 points}. The cell averages and point values at $t=1$ are displayed in Figure \ref{fig:KT}. We see that the solutions on this coarse mesh are consistent to \blue{and less oscillatory than}  what was obtained in our previous works, see e.g. in \cite{Abgrall2025a,Abgrall2024a}. 

\subsubsection*{Example 5---Double Mach Reflection Problem}
In the fifth example, the double Mach reflection (DMR) problem is used to test the robustness of the BP OE \pampa scheme.  The mesh contains 84,146 elements and  167,385 points. \red{Initially, a rightward shock with a Mach number of 10 is located at $x=-0.05$, and the initial condition is 
\begin{equation*}
  (\rho,u,v,p)=\left\{\begin{aligned}
  &(1.4,0,0,1),&&{x>-0.05},\\
  &(8,8.25,0,116.5), &&{x\leq-0.05}.
  \end{aligned}\right.
\end{equation*}}
Supersonic inflow boundary condition is specified on the left boundary, while free-stream outflow boundary condition is applied on the right boundary. At the upper boundary, the exact solution of
the Mach 10 moving oblique shock is imposed. Wall  conditions are imposed on the lower boundary.
 
\begin{figure}[ht!]
\centerline{\subfigure[Internal DoF $\xbar\rho_K$]{\includegraphics[trim=1.8cm 2.5cm 0.1cm 3.0cm,clip,width=5.5cm]{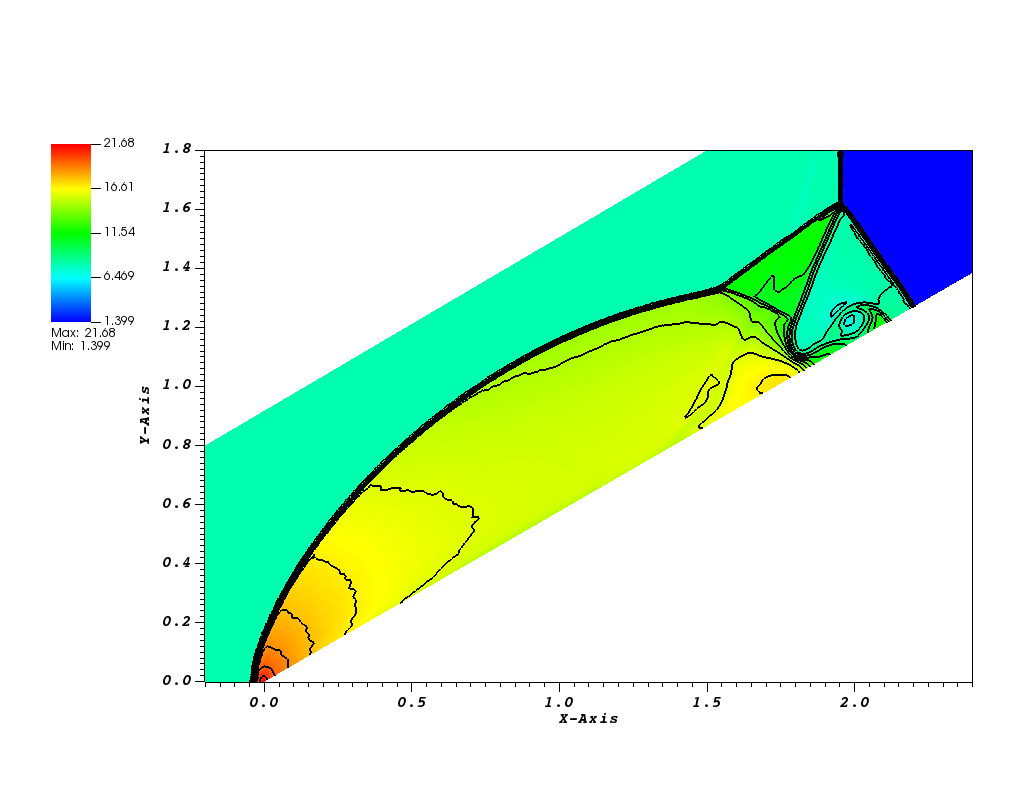}}\hspace*{0.05cm}
	\subfigure[Boundary DoFs $\rho_\sigma$]{\includegraphics[trim=1.8cm 2.5cm 0.1cm 3.0cm,clip,width=5.5cm]{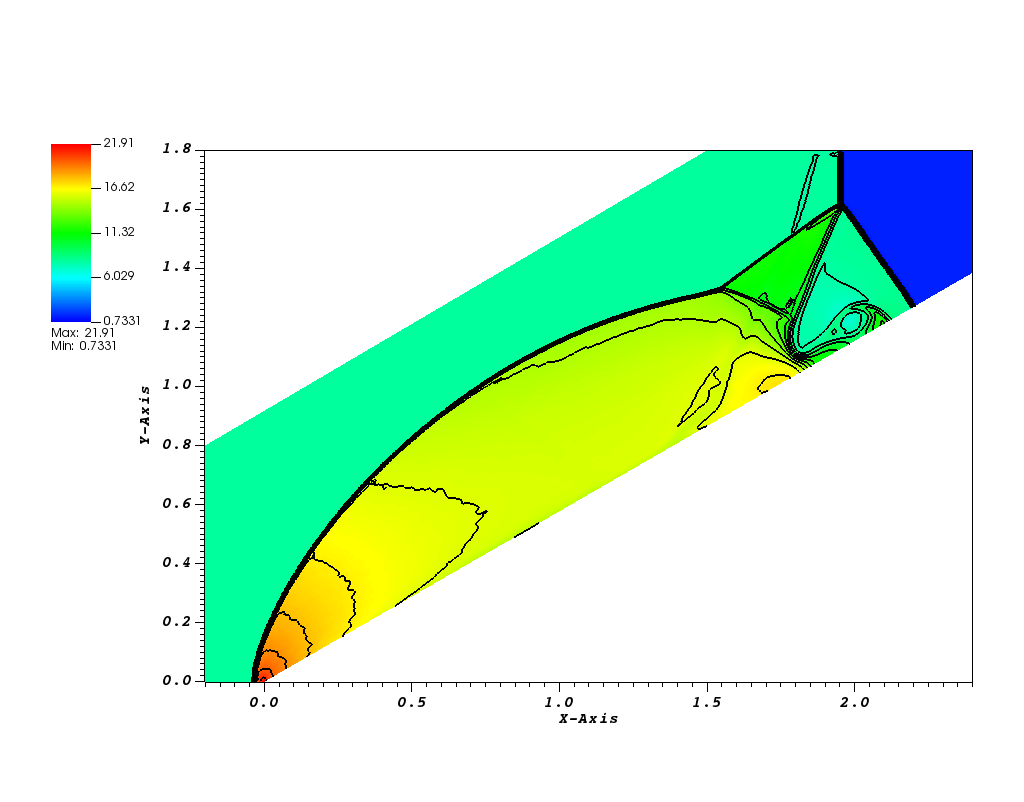}}}
\vskip1pt
\centerline{\subfigure[Blending parameter $\eta_e$]{\includegraphics[trim=1.8cm 2.5cm 0.1cm 3.0cm,clip,width=5.5cm]{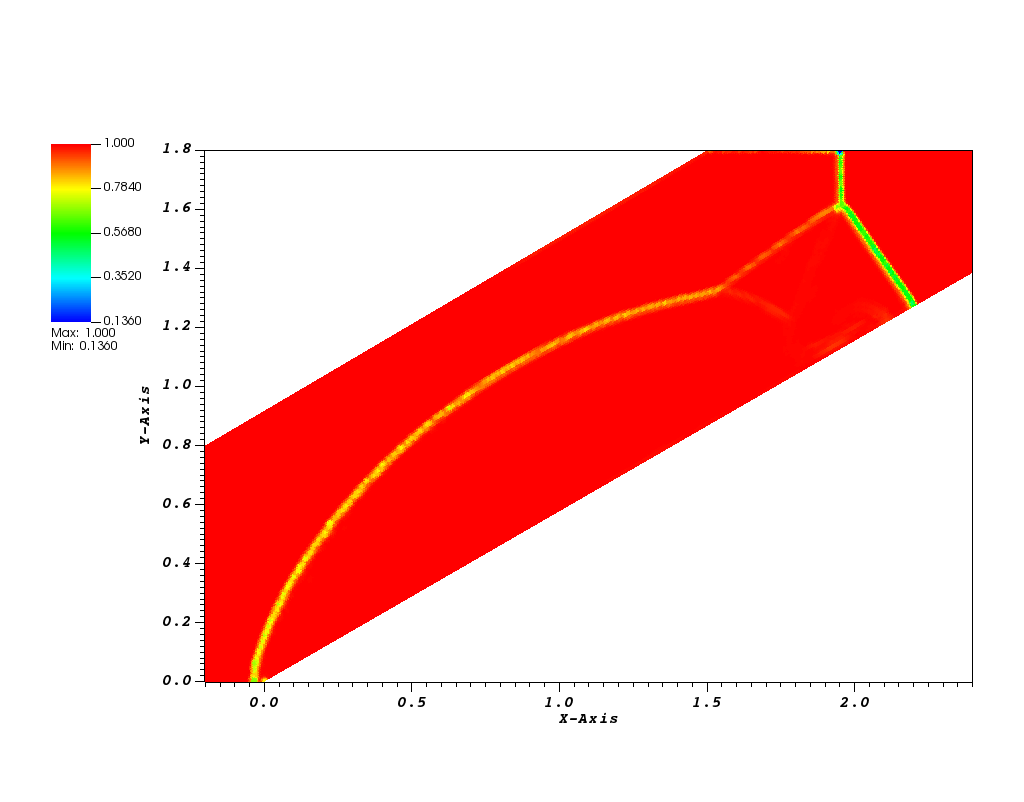}}\hspace*{0.05cm}
	\subfigure[Blending parameter $\eta_\sigma$]{\includegraphics[trim=1.8cm 2.5cm 0.1cm 3.0cm,clip,width=5.5cm]{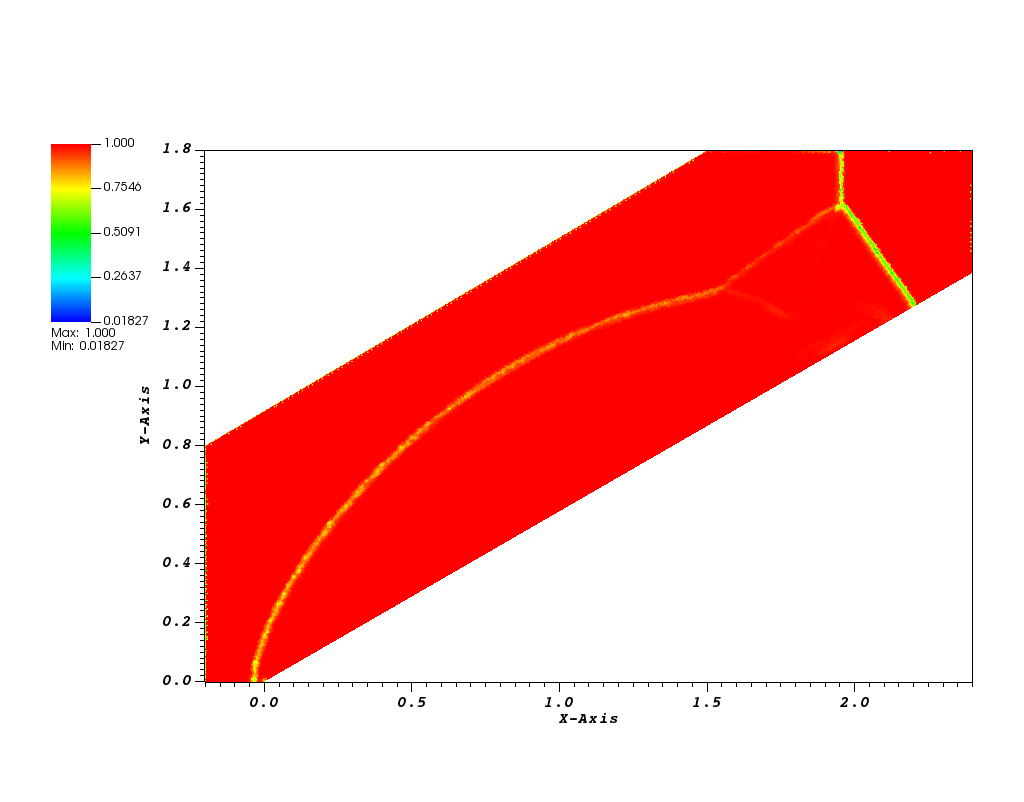}}}
\caption{\sf Example 5: DMR problem. Numerical solutions and BP OE blending parameters at $t=0.2$. \label{fig: DMR}}
\end{figure}

Figures \ref{fig: DMR}--(a) and (b) display the cell averages and point values of density, respectively. Compared with the results reported in \cite{Abgrall2025a}, the BP OE \pampa method produces noticeably fewer oscillations near shock waves, while remains more accurate profile around the double Mach stem. Furthermore, Figures \ref{fig: DMR}--(c) and (d) show the convex blending parameters for numerical fluxes and residuals. As expected, the blending between the high- and low-order schemes occurs primarily near discontinuities. Overall, these results demonstrate that the BP OE \pampa method effectively captures complex flow structures.

\subsubsection*{Example 6---Shock--Diffraction Problem}
In the sixth example, we simulate the diffraction of a shock wave at a convex corner with a $90^\circ$ angle. The computational domain is a union
of $[0,1]\times[-1,0]$ and $[-0.5,1]\times[0,1]$.  Initially, a pure right-going shock of Mach number 2.4, parallel to the $x$-axis, positioned at $x=-0.05$ and $0\leq y \leq 1$, is moving into undisturbed air ahead of the shock. The undisturbed air \red{(pre-shock)} has the state as
\begin{equation*}
  (\rho,\mbv,p)=(1,\mathbf 0, 1),
\end{equation*}
\red{while the post-shock state is defined according to the Rankine--Hugoniot relations.}
The boundary conditions are inflow at \red{$x=-0.5$}, $0\leq y\leq 1$, reflective at the walls \red{$-0.5\leq x\leq 0$}, $y=0$ and $x=0$, $-1\leq y\leq0$, and outflow on the remaining boundaries. 

\begin{figure}[ht!]
    \begin{center}
{\subfigure[Internal DoF $\xbar\bbu_K$]
{\includegraphics[width=0.45\textwidth]{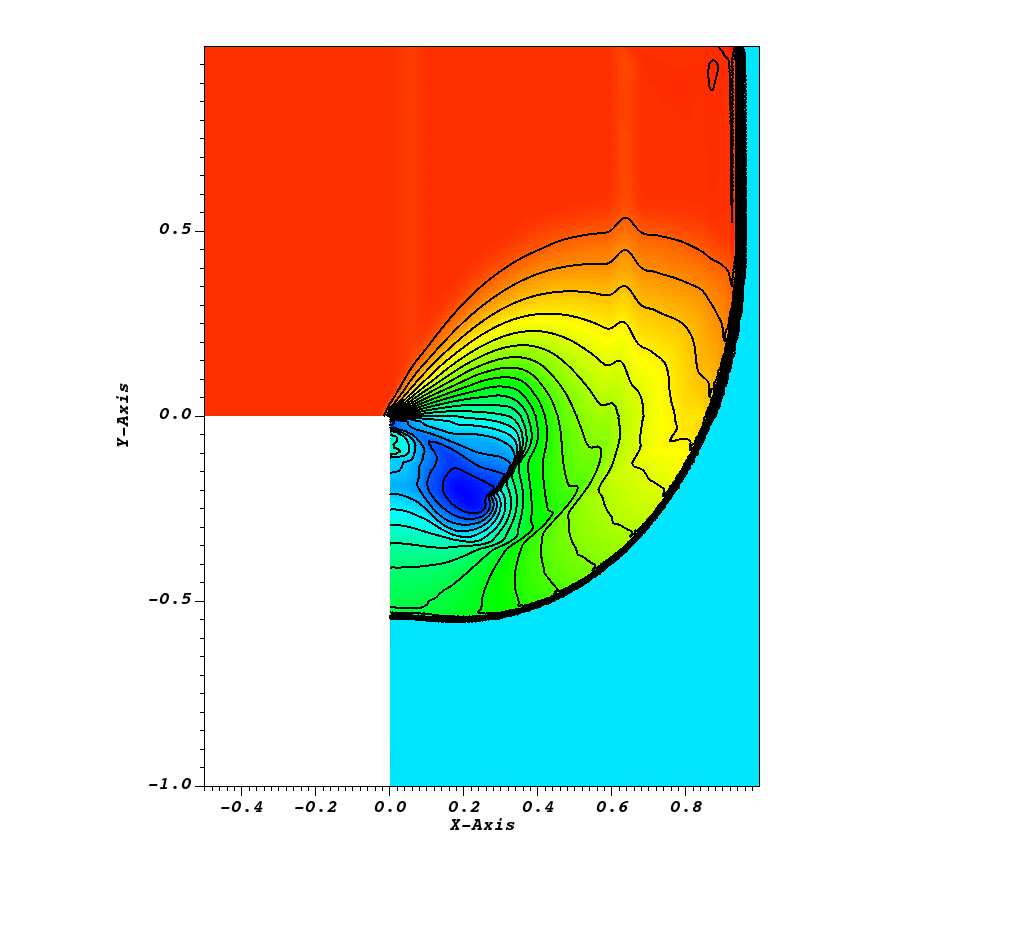}}}
\hspace*{0.05cm}
	\subfigure[Boundary DoFs $\bbu_\sigma$]{\includegraphics[width=0.45\textwidth]{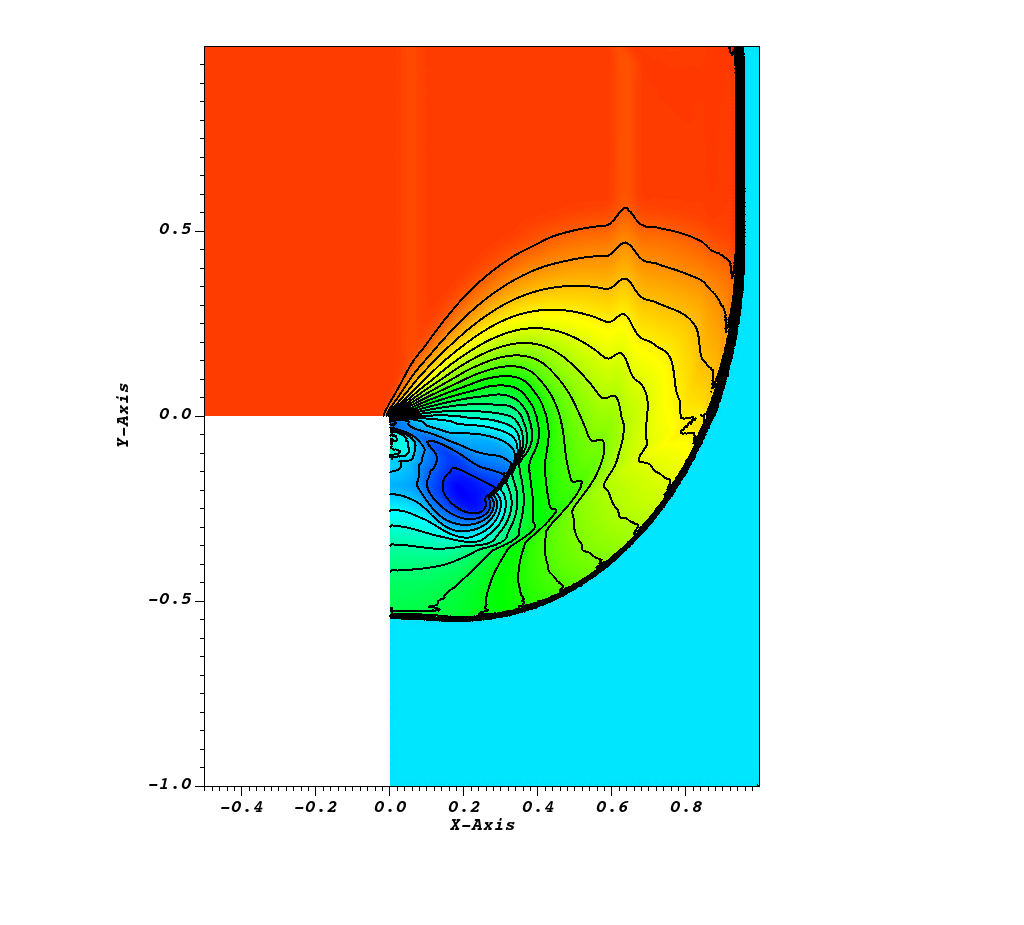}}
    {\subfigure[Sensor]{\includegraphics[width=0.45\textwidth]{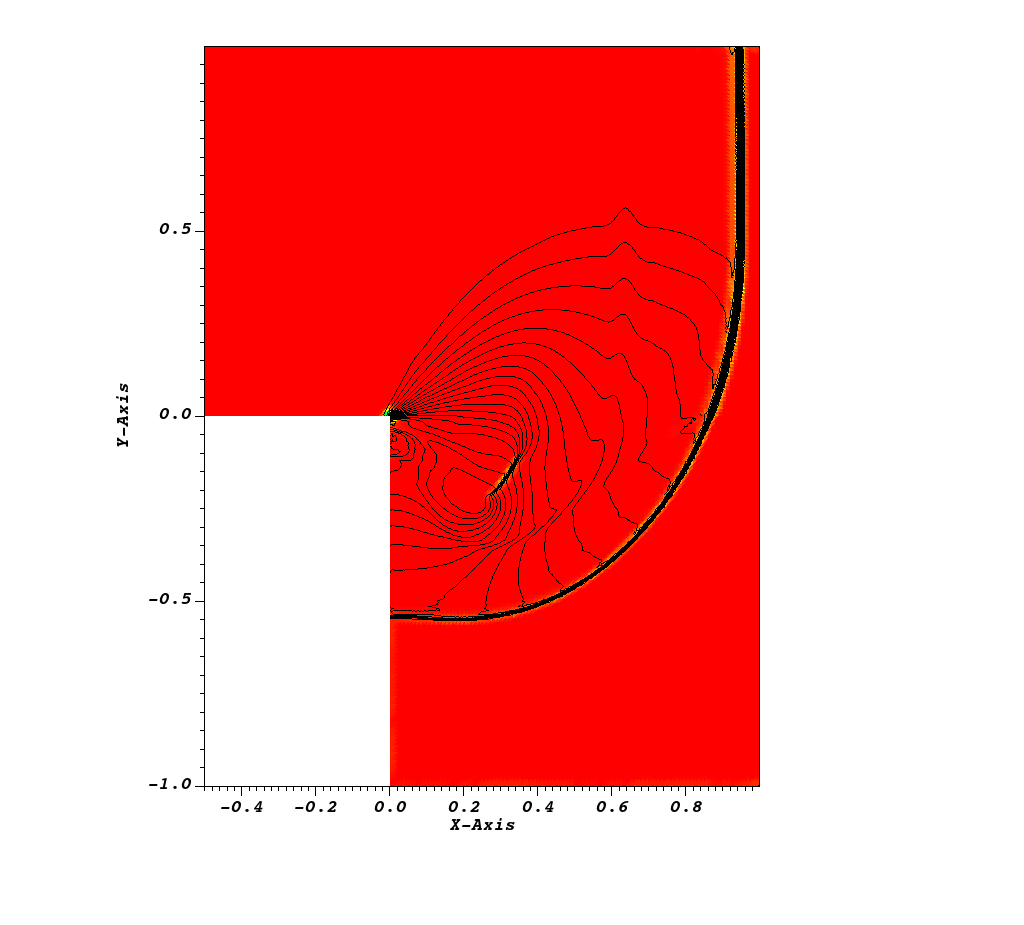}}}
    \end{center}
\caption{\sf Example 6: Shock--Diffraction problem. Numerical solutions (density) at $t=0.35$. \label{fig: Shock_Diffraction1}}
\end{figure}

Figure~\ref{fig: Shock_Diffraction1} displays the cell averages and point values of density computed on a mesh with 39,728 elements and 78,897 points at time $t=0.35$. With the BP OE blending parameter enabled, the scheme successfully captures the physical solution and the complex wave patterns in the cell averages. 

{\color{blue}
\subsubsection*{Example 7---Shock--Vortex Problem}
In the seventh example, we consider the interaction between an isentropic vortex and a Mach 1.1 shock, initially perpendicular to the $x$-axis and located at $x=0.5$. The initial left state of the shock is given by $(\rho,\bbv,p)=(1,1.1\sqrt{\gamma},0,1)$, while the right state is determined from the Rankine--Hugoniot conditions. An isentropic vortex, centered at $\bbx_c=(0.25,0.5)$, is superimposed onto the mean flow on the left of the shock. The perturbations in velocity, temperature ($T=p/\rho$), and entropy $(S=\ln(p/\rho^\gamma))$ due to the vortex are defined as
\begin{equation*}
    \delta\bbv=\frac{\varepsilon}{r_c}e^{\alpha(1-\eta^2)}(\bar y,-\bar x),\quad \delta T=-\frac{(\gamma-1)\varepsilon^2}{4\alpha\gamma}e^{2\alpha(1-\eta^2)},\quad \delta S=0,
\end{equation*}
where $r_c=0.05$, $\varepsilon=0.3$, $\alpha=0.204$. Here, $\eta=\frac{r}{r_c}$ with $r^2=\|\bar \bbx\|^2$ and $\bar\bbx=\bbx-\bbx_c$. The parameters $r_c$, $\varepsilon$, and $\alpha$ represent the critical radius, vortex strength, and decay rate, respectively. The computational domain is $[0,2]\times[0,1]$ and the mesh is generated by dividing each cell of a Cartesian $200\times100$ grid along its diagonal. Reflective boundary conditions are used at the upper and lower boundaries, while the left and right boundaries employ inflow and outflow conditions, respectively. We apply the BP OE \pampa-\dg scheme to run the simulation up to final time $t=0.8$. The contour plots of pressure at six different time snapshots are shown in Figure \ref{fig: Shock_Vortex}. The results are comparable to those reported in \cite{Lu2021,Peng2025}. The shock-vortex interaction dynamics and the intricate wave structures are well captured by the proposed BP OE \pampa-DG scheme.
\begin{figure}[ht!]
\centerline{\subfigure[$t=0.068$]{\includegraphics[trim=4.2cm 2.5cm 1.5cm 1.0cm,clip,width=4.2cm]{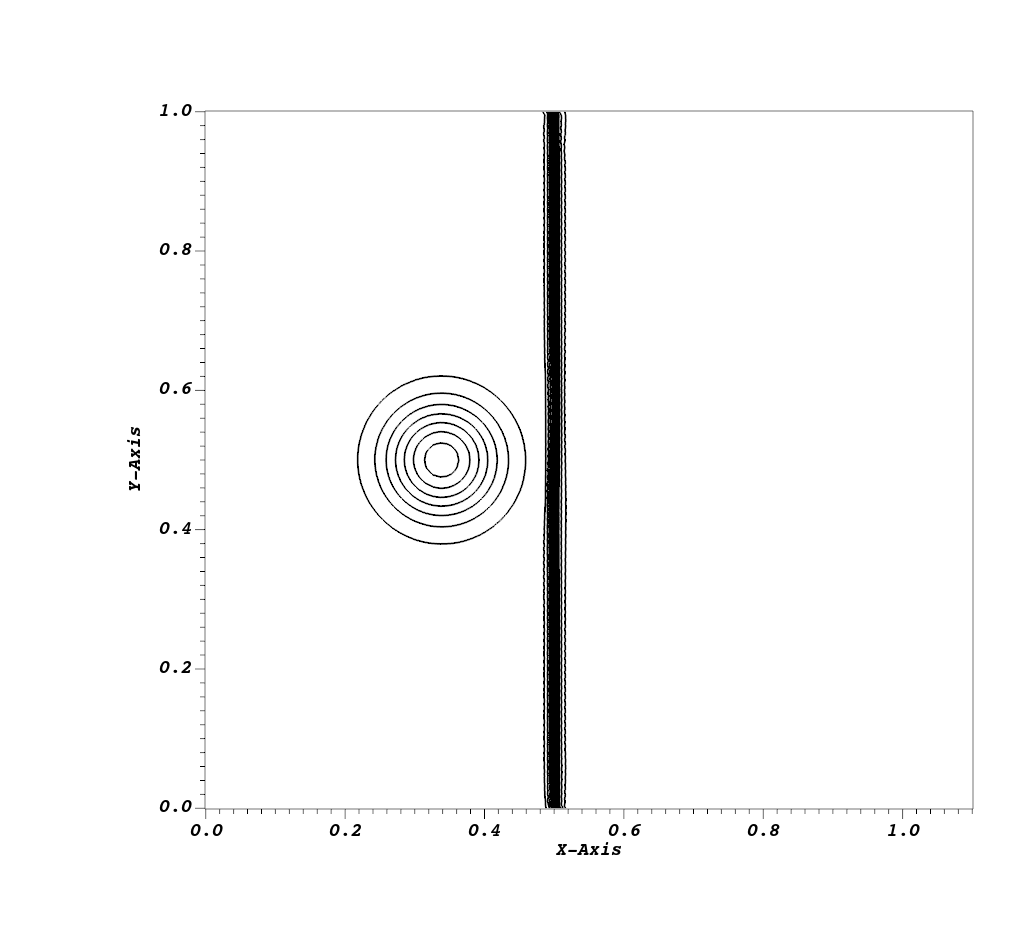}}\hspace*{0.05cm}
	\subfigure[$t=0.203$]{\includegraphics[trim=4.2cm 2.5cm 1.5cm 1.0cm,clip,width=4.2cm]{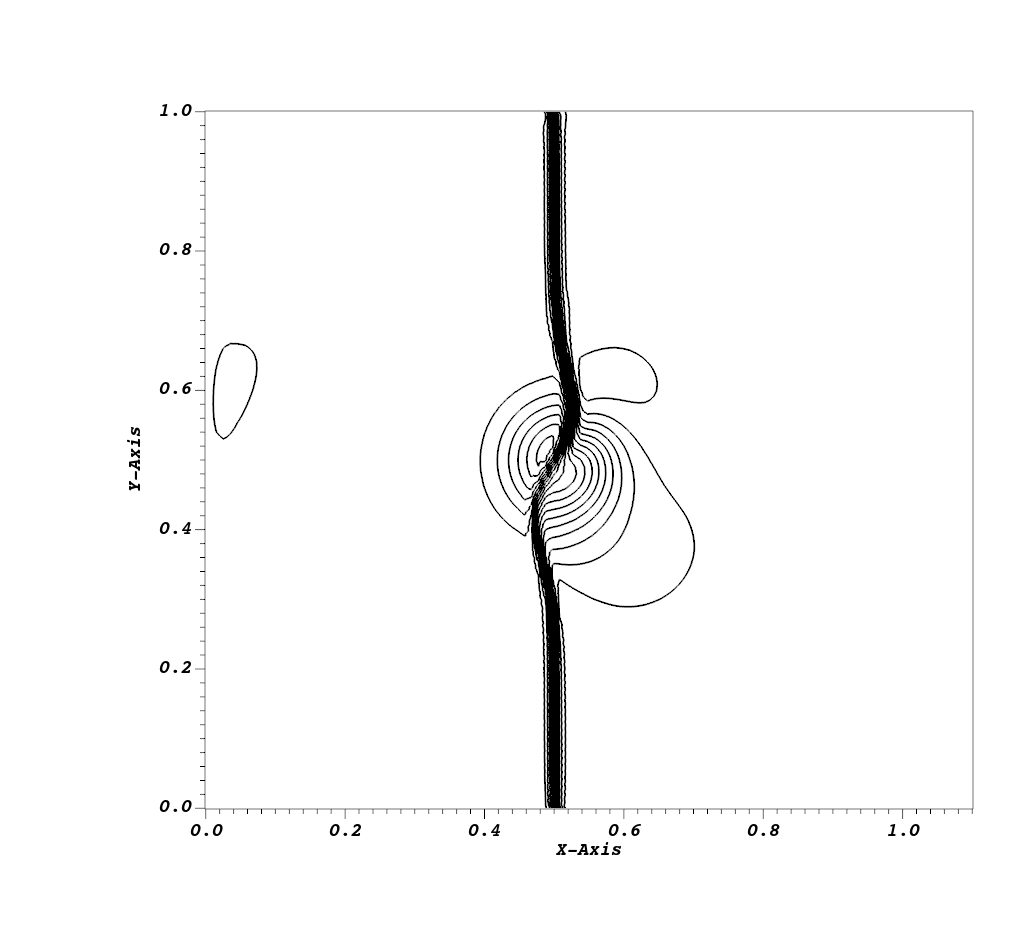}}\hspace*{0.05cm}
	\subfigure[$t=0.330$]{\includegraphics[trim=4.2cm 2.5cm 1.5cm 1.0cm,clip,width=4.2cm]{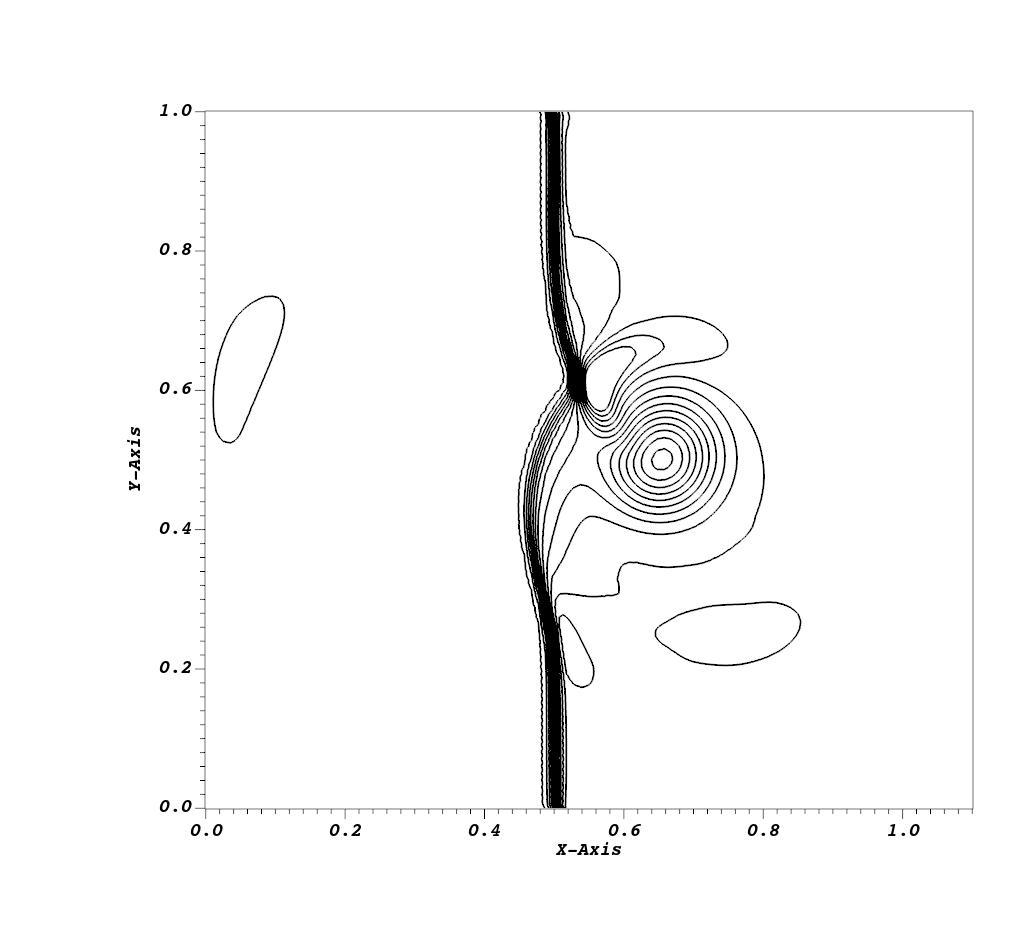}}}
    \vskip1pt
    \centerline{\subfigure[$t=0.529$]{\includegraphics[trim=4.2cm 2.5cm 1.5cm 1.0cm,clip,width=4.2cm]{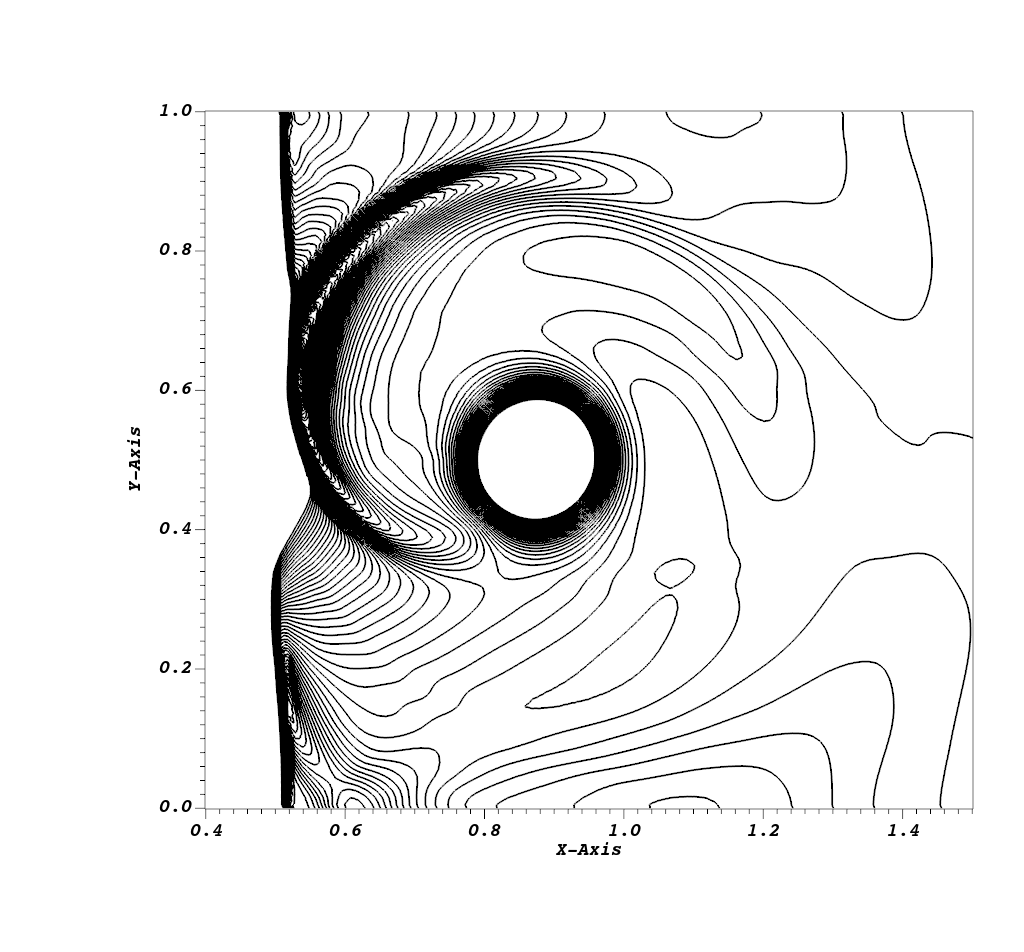}}\hspace*{0.05cm}
	\subfigure[$t=0.662$]{\includegraphics[trim=4.2cm 2.5cm 1.5cm 1.0cm,clip,width=4.2cm]{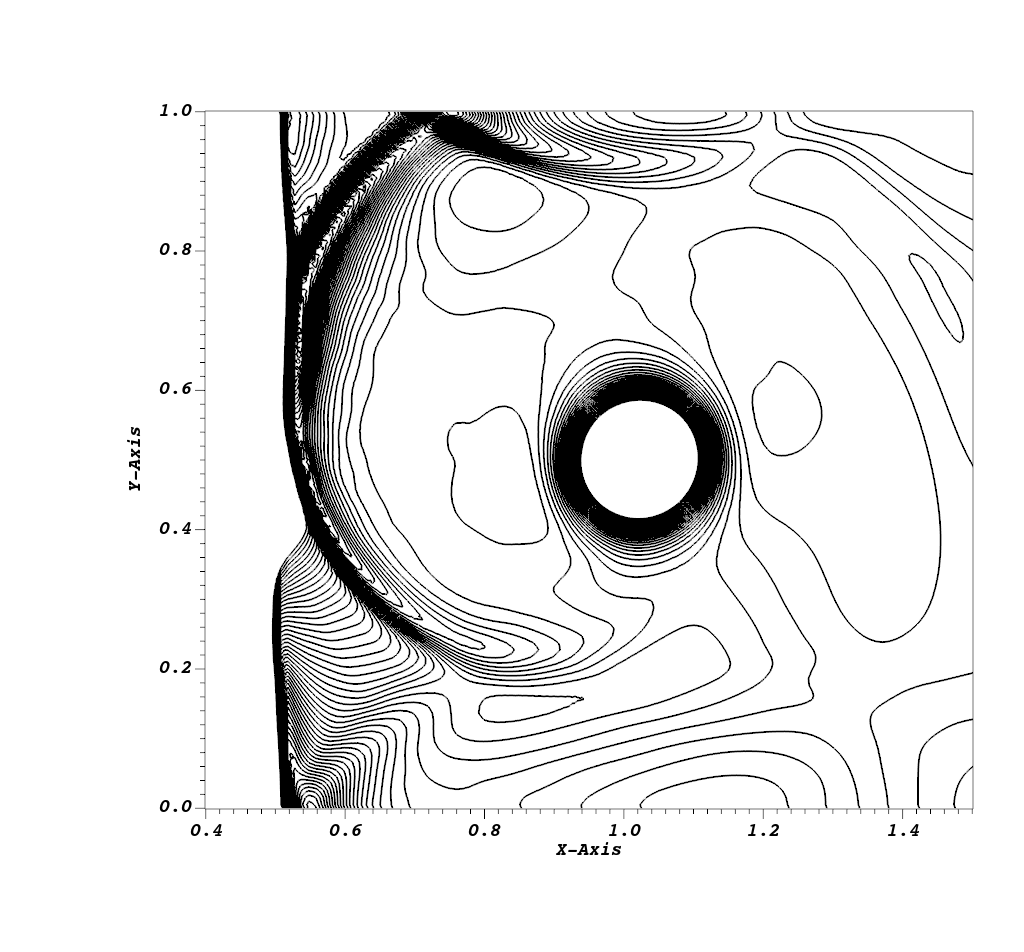}}\hspace*{0.05cm}
	\subfigure[$t=0.8$]{\includegraphics[trim=4.2cm 2.5cm 1.5cm 1.0cm,clip,width=4.2cm]{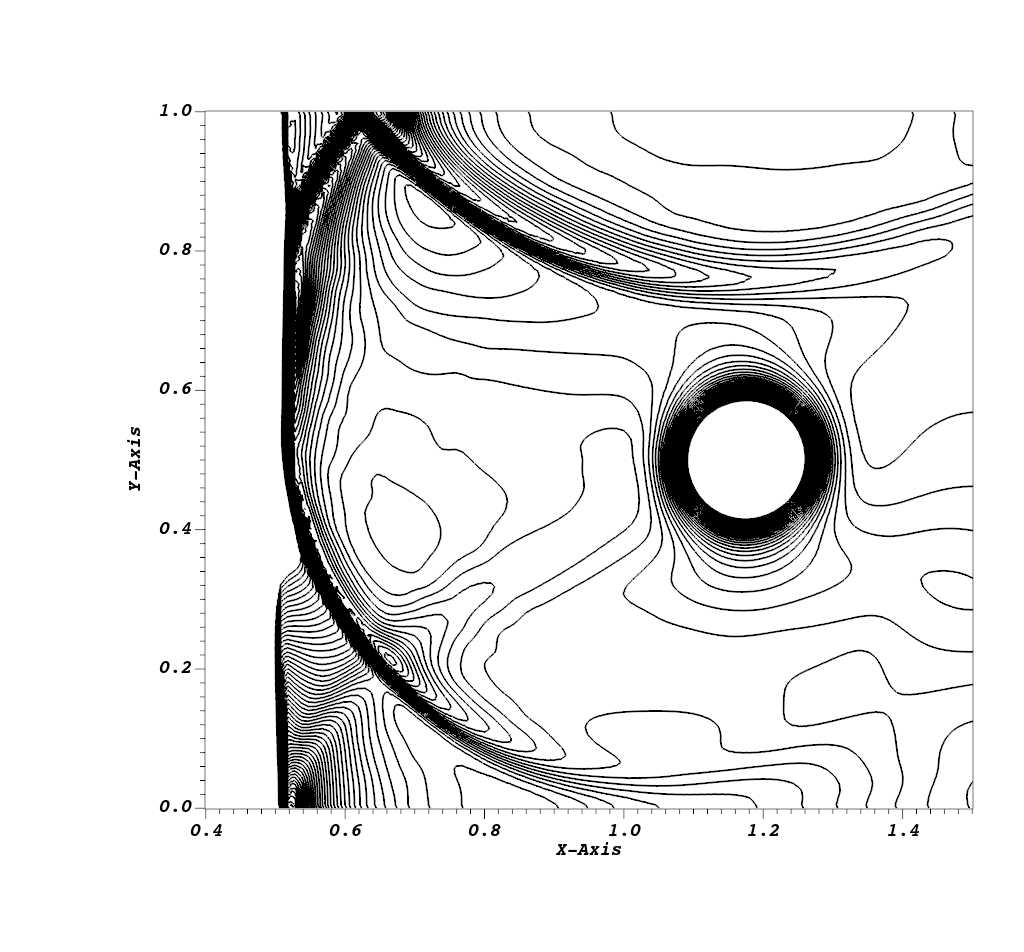}}}
\caption{\sf Example 7: Shock--Vortex problem. Contour plots of pressure in cell averages. Top row: 30 contour lines in [0.69, 1.3]; Bottom row: 90 contour lines in [1.19, 1.36]. \label{fig: Shock_Vortex}}
\end{figure}

\subsubsection*{Example 8---Sedov Blast Wave Problem}
In the eighth example, we consider the Sedov problem \cite{sedov}, which is characterized by extremely low density and very strong shock waves. The computational domain is a circle centered at $(0,0)$ of radius $1.5$. We impose the outflow boundary conditions on all sides. The initial conditions are prescribed as follows:
\begin{itemize}
    \item For all cells, we set $(\rho,\bbv,E)(x,0)=(1,\mathbf 0,e_{\min})$ with $e_{\min}=10^{-12}$;
    \item Except for the  cells containing the origin $(0,0)$, we set $(\rho,\bbv,E)(x,0)=(1,\mathbf 0, e_{\max})$ where $e_{\max}=0.979264/\vert K\vert$, and $\vert K\vert$ denotes the area of the cell.
\end{itemize}
These initial conditions are chosen so that the shock strength is almost infinite, and the initial condition corresponds to a Dirac mass or energy at $\bbx=\mathbf 0$ in the limit of mesh refinement. This setup generates a strong radially symmetric blast wave, providing a stringent test for the robustness and positivity-preserving capability of the numerical scheme. We display the obtained numerical results at $t=1$ in Figure \ref{fig: sedov}. 
\begin{figure}[ht!]
\centerline{\subfigure[Internal DoF $\xbar\rho_K$]{\includegraphics[trim=1.2cm 2.5cm 1.5cm 1.0cm,clip,width=4.2cm]{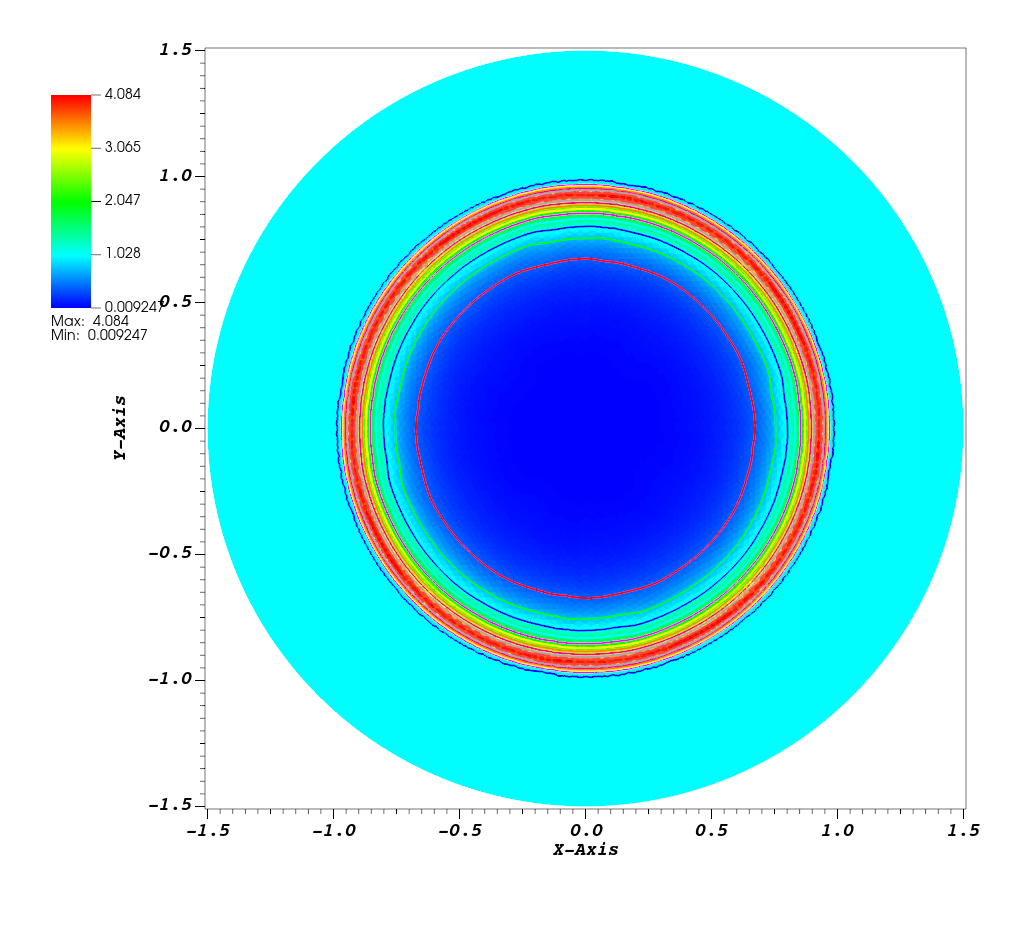}}\hspace*{0.05cm}
	\subfigure[Boundary DoFs $\rho_\sigma$]{\includegraphics[trim=1.2cm 2.5cm 1.5cm 1.0cm,clip,width=4.2cm]{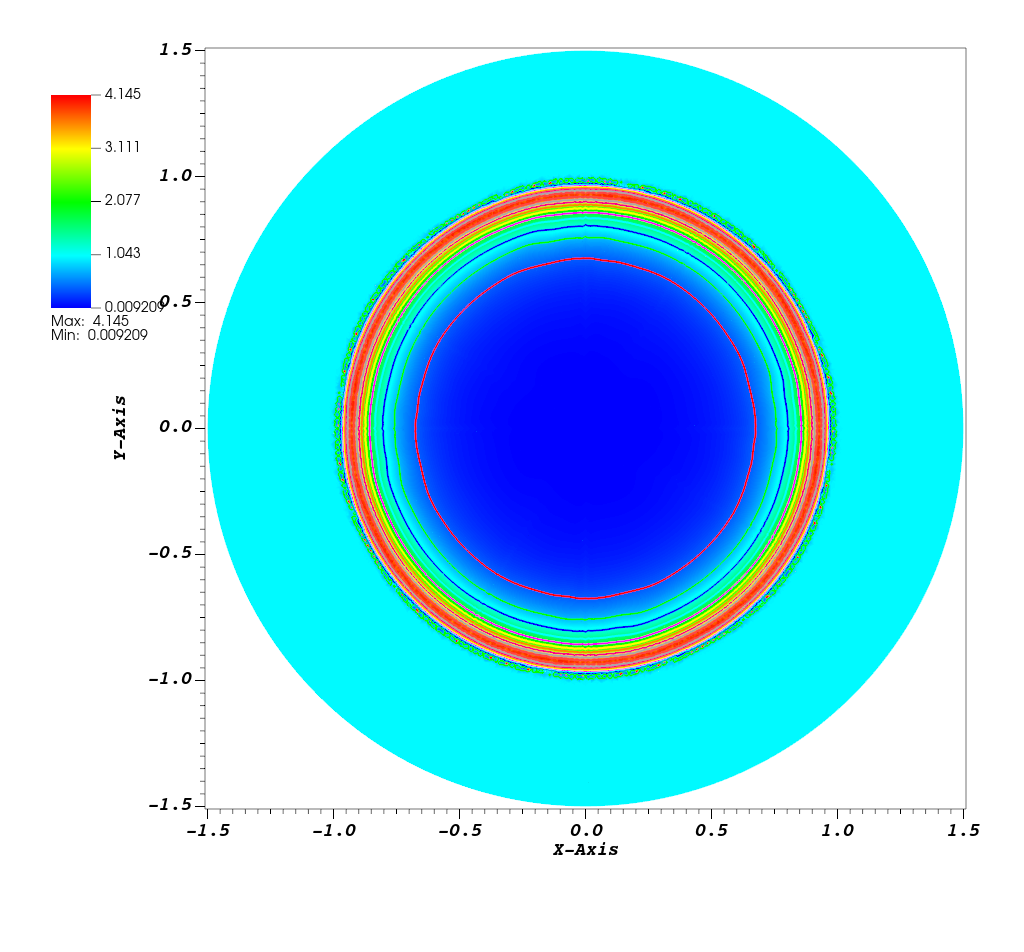}}\hspace*{0.05cm}
	\subfigure[Cut line along $y=x$]{\includegraphics[trim=1.2cm 2.5cm 1.5cm 1.0cm,clip,width=4.2cm]{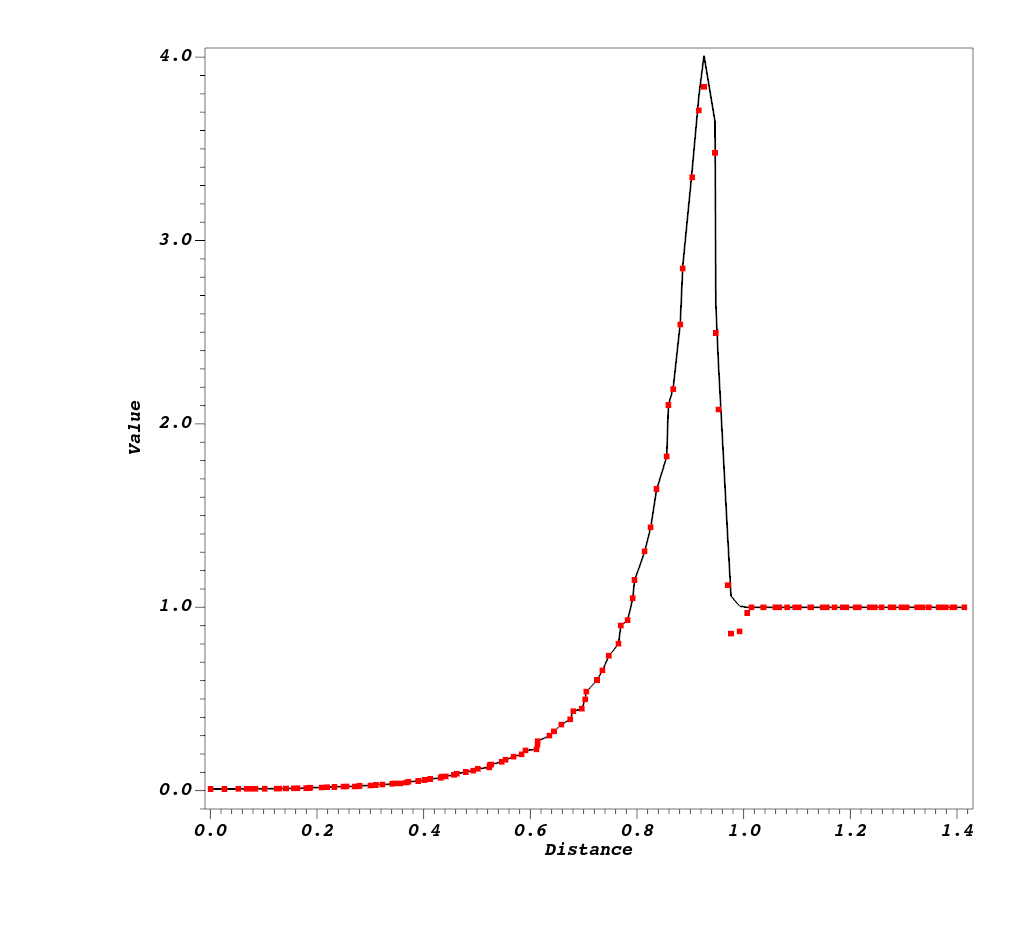}}}
\caption{\sf Example 8: Sedov blast wave problem. Density at time $t=1$. \label{fig: sedov}}
\end{figure}

\subsubsection*{Example 9---Mach 2000 Astrophysical Jet Problem}
In the ninth example, we consider the challenging astrophysical jet problem of Mach number 2000, following \cite{Zhang2010a}. A computational domain $[0,1]\times[-0.25,0.25]$ is considered with a fluid at conditions $(\rho,\bbv,p)=(0.5,\mathbf 0,0.4127)$. A high-speed jet state $(\rho,\bbv,p)=(5,800,0,0.4127)$ is injected into the domain from the left boundary within the range $y=-0.05$ to $0.05$. The remaining boundaries are outflow conditions. The results at $t=0.001$ are displayed in Figure \ref{astrojet}, as well as the flag on elements (the one on points has a similar behavior), and a zoom of the mesh.
\begin{figure}[ht!]
\centerline{
    \subfigure[Boundary DoFs $\log(\rho_\sigma)$]{\includegraphics[width=0.45\textwidth]{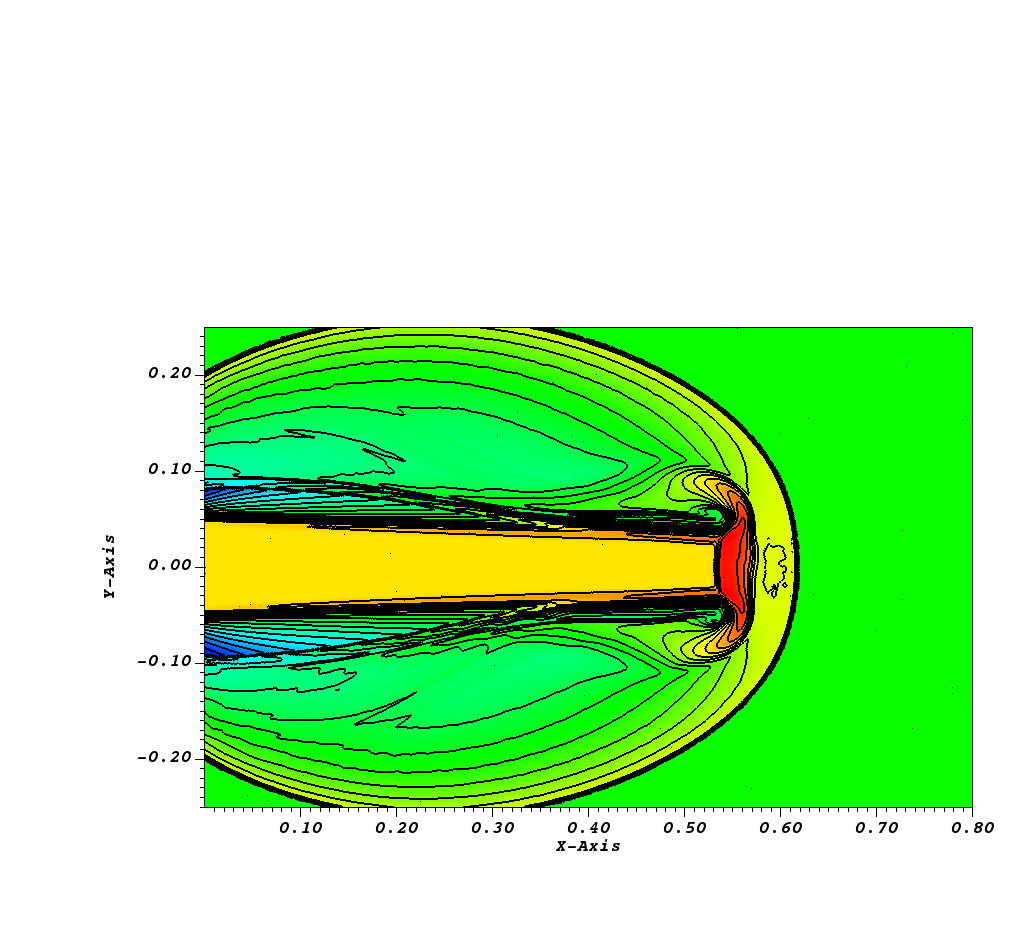}}
    \subfigure[Internal DoFs $\log(\xbar{\rho}_K)$]{\includegraphics[width=0.45\textwidth]{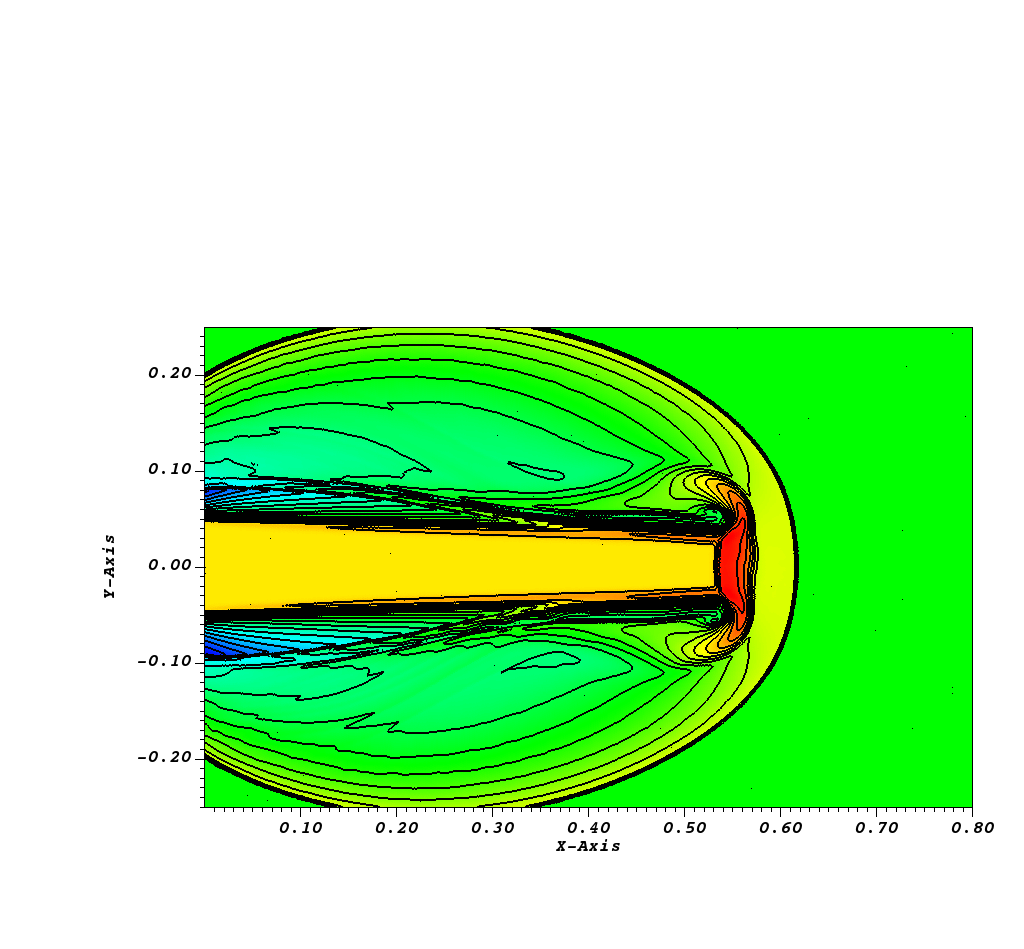}}}
    \vskip1pt
    \centerline{
    \subfigure[Flag Elements]{\includegraphics[width=0.45\textwidth]{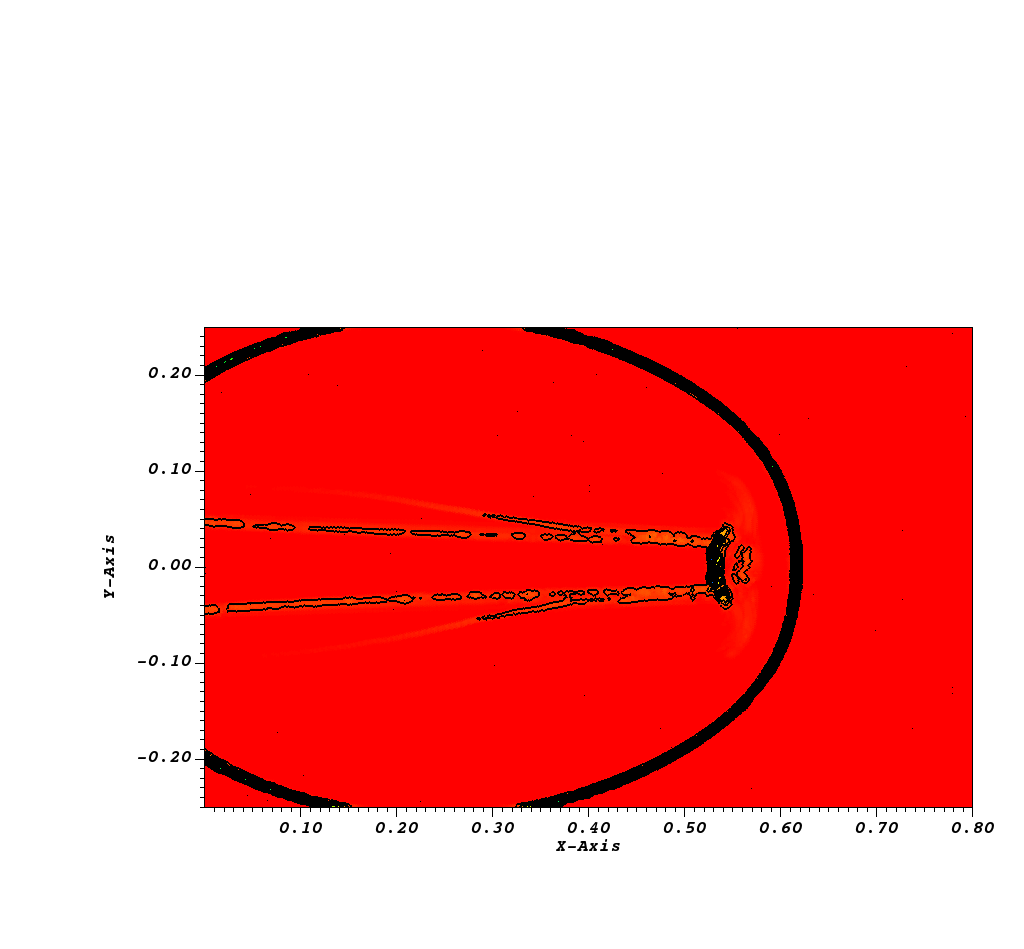}}
    \subfigure[Boundary elements and mesh]{\includegraphics[width=0.45\textwidth]{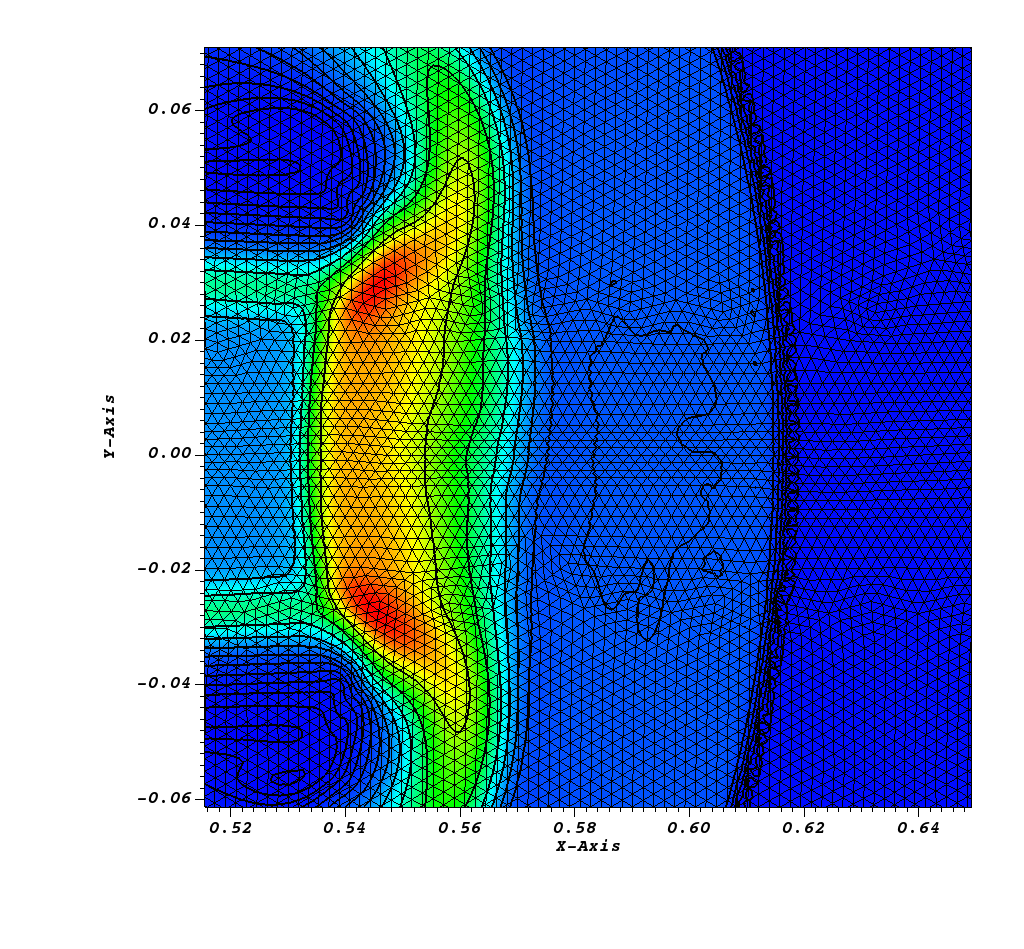}}}
\caption{Example 9: Astro jet problems. The ranges are $\rho_\sigma\in [0.007669, 32.14]$ and $\xbar{\rho}\in [0.006046,32.35]$. 20 isolines are used in logarithm coordinates. \label{astrojet}}
\end{figure}
In this mesh we have $465,889$ point values DoFs and $232,320$ internal DoFs. This corresponds to approximately to a $350\times 640$ Cartesian mesh.

\begin{remark}[About the boundary conditions] In that case, the top, bottom and right boundary conditions are outflow ones. The more complicated one is on the left boundary from which the jet enters the domain. The natural condition is to set the normal derivatives equal to 0. This is very easy to do for a Cartesian mesh. In the case of a triangular type mesh, which is not obtained by cutting into two triangles quadrangles, it is a bit more complex. We have not been able to find a good way to set up homogeneous Neumann conditions because the problem is written with first order differential operators. After many tries, the best solution seems to enlarge the domain on the left: is is now $[-0.1,0.8]\times [-0.25,0.25]$. Then we impose strongly the boundary conditions on the left ($x=-0.1$) and post-process the solution in the buffer zone $[-0.1,0]\times [-0.25,0.25]$ by blending the physical variables obtained at each Runge Kutta cycle with the initial solution. We have chosen, after several tries the following blending parameter: at the point $(x,y)$ we set
$$\theta(x,y)=\left \{\begin{array}{ll}1 &\text{if }x\geq -0.05\\
\big (\frac{\max(\vert x+0.05\vert,\vert y+0.05\vert)}{0.1}\big )^5 & \text{else}
\end{array}\right .$$
so that we keep intact the solution for $x\geq -0.05$.
\end{remark}

\subsubsection*{Example 10---Flow Past a Circular Cylinder} 
We conclude the numerical validations with the benchmark problem, which involves a rightward parallel flow with a Mach number of 25 passing a circular cylinder within the domain $[-3,0]\times[-6,6]$. The circular cylinder, with a unit radius, is centered at the origin on the $x$-$y$ plane. The initial conditions are $(\rho,\bbv,p)=(1.4,25,0,1)$. Inflow boundary condition is applied along the left boundary, reflective boundary conditions are imposed on the top and bottom boundary. Wall conditions are imposed on the circle.

This problem is steady, so that, in principle, the solution should be independent of the CFL. Note that the time step appears in \eqref{eq:OE_param}, so we have used the following form of $\theta_K$:
$$
    \theta_K=\exp\Big(\frac{-1}{N_K}\sum_{e\subset\partial K}
    \sigma_{e,K}(\mbu_{\mathrm h})\Big).
$$
This amounts to setting 
$\tfrac{\alpha_e\dt_n}{\ell_{e,K}}\equiv 1$
and can be interpreted as follows: when one computes a steady problem, generally one uses a local time stepping method: this means that the time step is different for each cell/point, and is computed so that the local CFL number is uniform.  Noticing that $\frac{\alpha_e\dt_n}{\ell_{e,K}}$ is a local CFL number, we proceed in \eqref{eq:OE_param} as if the CFL number were equal to 1.

In Figure \ref{blunt}, we show the density, the mesh, the sensor on the elements (it is similar for the points), and the density residual. The mesh is rather coarse: we have $11,541$ nodes and $5,628$ elements. In total, we have $22,796$ boundary DoFs and $5,628$ internal DoFs.

\begin{figure}[ht!]
\centerline{\subfigure[Density with 20 iso-lines]{\includegraphics[trim=1.8cm 2.5cm 0.1cm 3.0cm,clip,width=5.5cm]{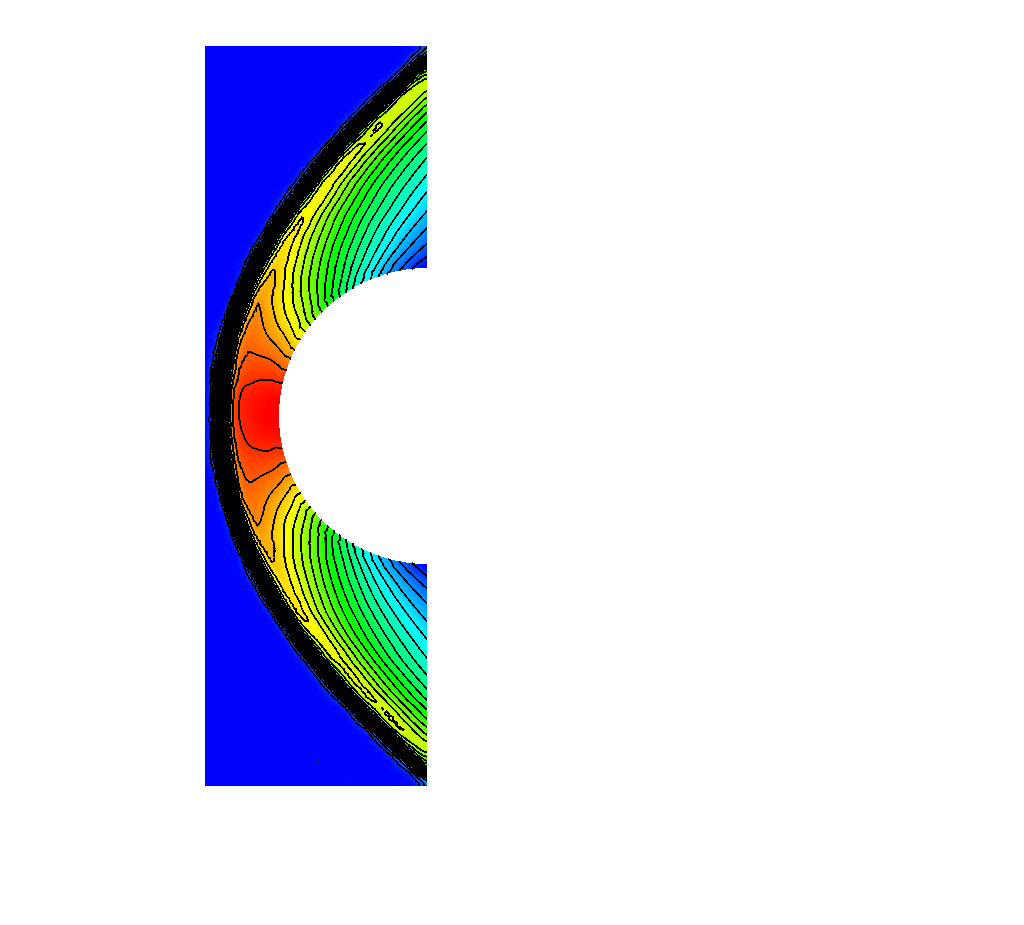}}\hspace*{0.05cm}
	\subfigure[Density and sensor on the elements]{\includegraphics[trim=1.8cm 2.5cm 0.1cm 3.0cm,clip,width=5.5cm]{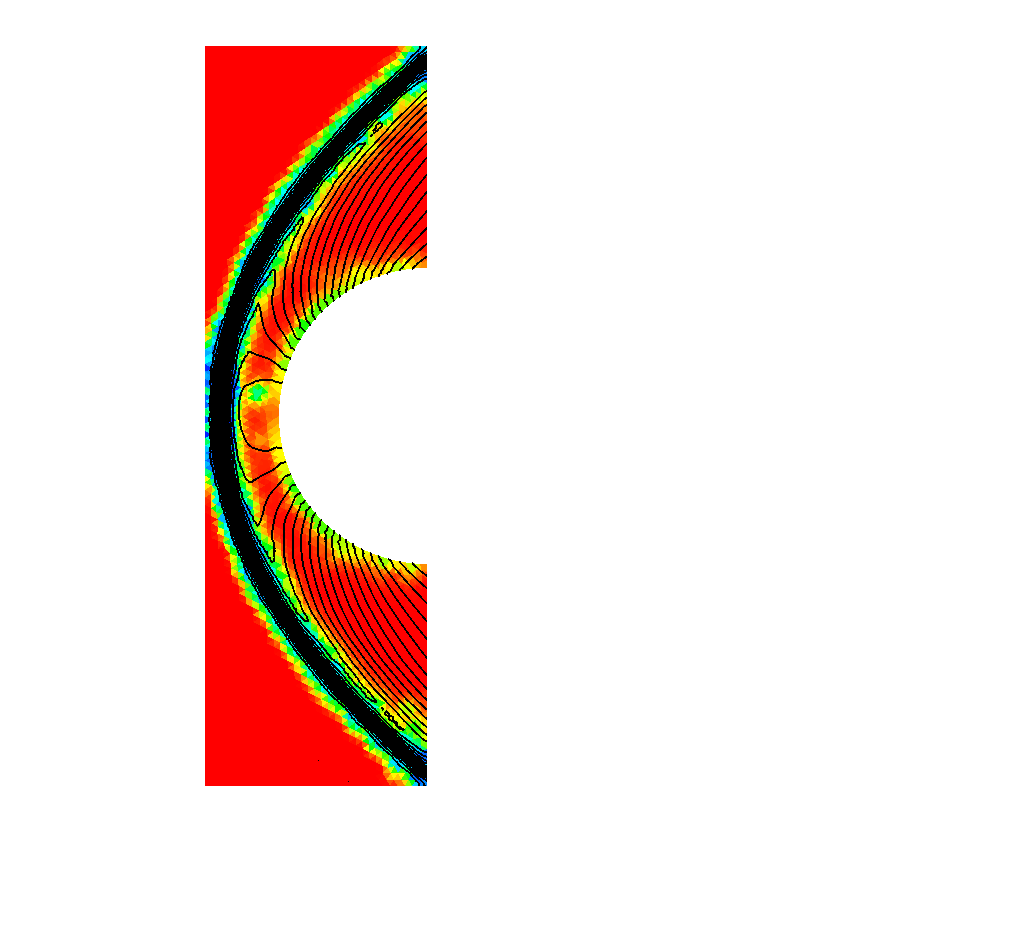}}}
\vskip1pt
\centerline{\subfigure[Density residuals]{\includegraphics[trim=1.8cm 2.5cm 0.1cm 3.0cm,clip,width=5.5cm]{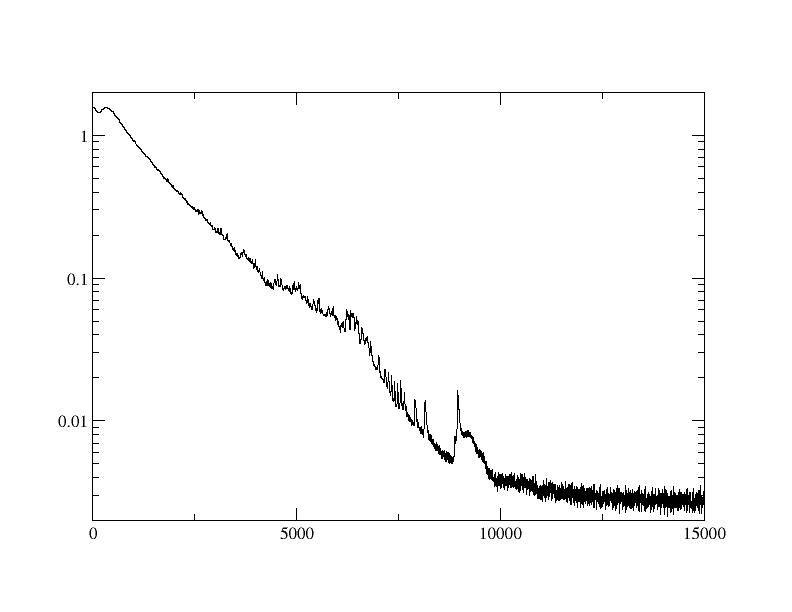}}\hspace*{0.05cm}
	\subfigure[zoom of the density and the mesh]{\includegraphics[trim=1.8cm 2.5cm 0.1cm 3.0cm,clip,width=5.5cm]{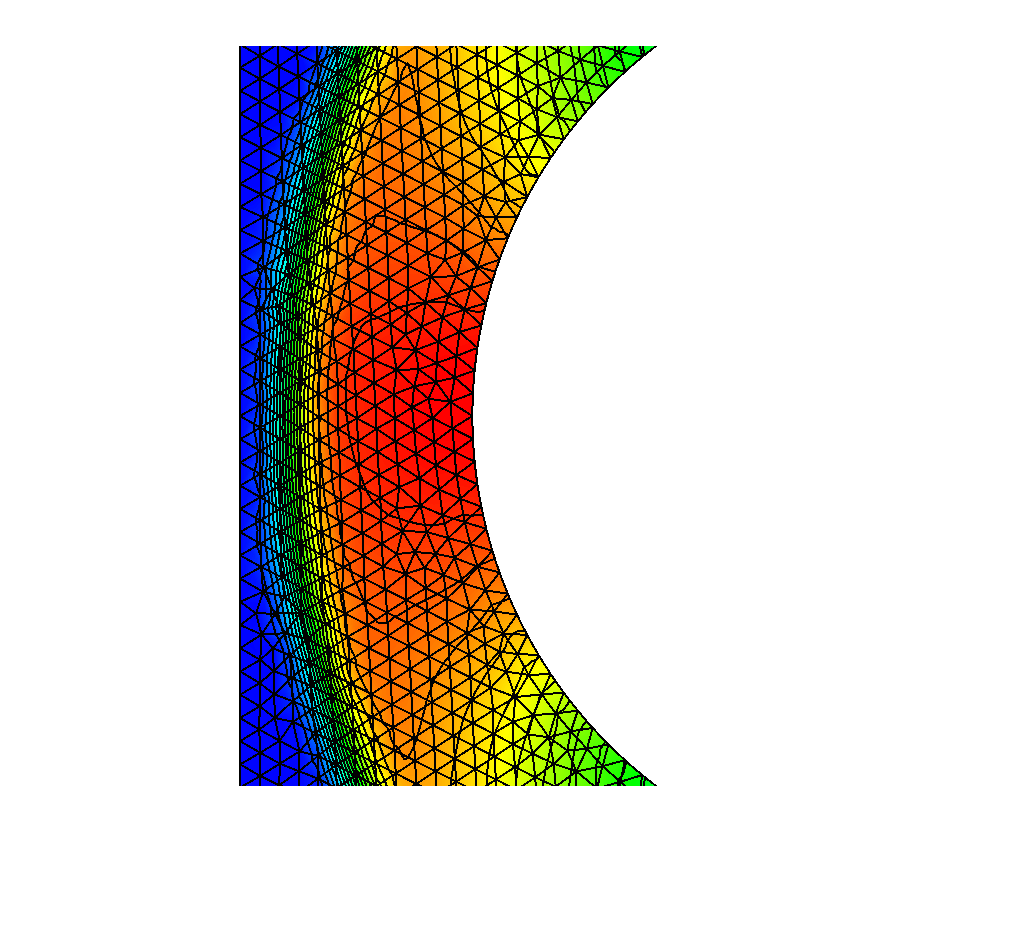}}}
\caption{\sf Example 10: Flow past a circular cylinder problem. Numerical solutions. The range for density is $[1.365,8.649]$. \label{blunt}}
\end{figure}
We can see, on Figure \ref{blunt}-(c), that the residual stagnates after 10000 time steps: nothing has been do to have an efficient method for this kind of problem, so there is big surprise here. However, this method is again shown to be very robust. 
}

\section{Conclusions and Perspectives}
We have developed a compact, third-order accurate, bound-preserving, and non-oscillatory \pampa scheme for hyperbolic conservation laws, specifically for scalar conservation laws and the Euler equations of gas dynamics, using a \dg formulation on unstructured triangular meshes. By properly selecting the projector weights, the proposed method recovers the upwind \pampa scheme previously introduced in \cite{Abgrall2024a,Abgrall2026}.

To ensure the bound preservation of both internal and boundary DoFs, we employ a sufficient convex BP parameter originally derived in \cite{Abgrall2024a} to blend the high- and low-order fluxes and residuals. To further enhance robustness in the presence of strong shocks, we introduce a problem-independent convex oscillation eliminating parameter based on a rotationally invariant damping term from \cite{Ding2025}. By taking the minimum of the BP and OE parameters, we obtain a unified convex BP OE blending parameter. 
Extensive numerical experiments, including smooth and discontinuous test cases and a range of models from scalar conservation laws to Euler systems of gas dynamics, demonstrate the stability and accuracy of our proposed method.

\blue{Additionally, we have also provided the truncation error analysis to confirm the third-order accuracy as expected and observed from the numerical experiments. The study of other properties, such as entropy properties, error estimates, and superconvergence, similar to those of the \dg scheme,  have not yet been done. Of course these are the next things on our to-do list.}

Other directions in our future work will focus on alternative projector constructions, such as those based on arithmetic averaging, which lead to a central \pampa scheme. However, preliminary numerical evidence (not shown here) indicates that the central \pampa scheme may require additional stabilization. Moreover, we aim to extend the framework to higher-order methods via higher-moment DoFs and to general polygonal meshes, following ideas sketched in \cite{Abgrall2025d}.

\bibliographystyle{siamplain}
\bibliography{refer_list}

\appendix


\section{Truncation error analysis}\label{sec:Truncation_error}
We consider the problem \eqref{eq:1} in $\R^d$
\begin{subequations}
\label{appendix:eq:1}
\begin{equation}
\label{appendix:eq:1:1}
\dpar{\bbu}{t}+\text{div }{\bbf(\bbu)}=0, \qquad \bbx\in \R^d
\end{equation}
with the initial condition 
\begin{equation}
\label{appendix:eq:1:2}
\bbu(x,0)=\bbu_0(x), \qquad \bbx\in \R^d.
\end{equation}
Here $\bbu\in \mathcal{D}\subset \R^m$. 
\end{subequations}

We assume $d=2$ to fix ideas, the case $d=3$ would be done the same. We consider a triangulation of $\R^2$, made of triangles. The triangles are denoted generically by $K$. The vertices of the triangles are denoted by $\bba_j$, and the degrees of freedom (hence including these vertices) by $\sigma_l$. The centroid of $K$ is $\bbx_K$.

The scheme for $\xbar\bbu_K$ writes
\begin{subequations}
\label{appendix:scheme}
\begin{equation}
\label{appendix:scheme:1}
\vert K\vert   \; \dfrac{\mathrm d\xbar \bbu_{K}}{\mathrm dt}+  \oint_{\partial K}\bbf(\bbu_{\mathrm h})\cdot\bbn \; \mathrm d\ell  =0.
\end{equation}
For $\bbu_\sigma$, we assume a semi-discrete scheme of the following form:
\begin{equation}
\label{appendix:scheme:3}
 \dfrac{\mathrm d\bbu_\sigma}{\mathrm dt}+ \sum_{K,\sigma\in K}\omega_{K,\sigma} \Phi_{K,\sigma}(\bbu_{\mathrm h})=0,
\end{equation} 
\end{subequations}
where in $K$, we have
\begin{equation}
\label{appendix:zeschema}\vert K\vert\;
\begin{pmatrix} \Phi_{K,\sigma_1}(\bbu_{\mathrm h})\\ \vdots \\ \Phi_{K,\sigma_6}(\bbu_{\mathrm h})\\[0.5em] \frac{1}{\vert K\vert}\oint_{\partial K}\bbf(\bbu_{\mathrm h})\cdot\bbn \; \mathrm d\ell\end{pmatrix}
=\PP \begin{pmatrix} \vdots \\ -\int_{K}\nabla\varphi_{\sigma_l}\cdot\bbf(\bbu_{\mathrm h})\;\mathrm d\mbx+\oint_{\partial K}\varphi_{\sigma_l}\bbf(\bbu_{\mathrm h})\cdot \bbn \; \mathrm d\ell \\ \vdots \\-\int_{K}\nabla\xbar\varphi\cdot\bbf(\bbu_{\mathrm h})\;\mathrm d\mbx+\oint_{\partial K}\xbar \varphi\bbf(\bbu_{\mathrm h})\cdot \bbn \; \mathrm d\ell\end{pmatrix},
\end{equation}
where $\PP={\vert K\vert}M_K^{-1}$ as defined in \eqref{eq:Project}.
The matrices/scalars $\omega_{K,\sigma}$ satisfy
$$\sum_{K, \sigma\in K}\omega_{K,\sigma}=\Id \text{ and }\Vert \omega_{K,\sigma}\Vert \leq C$$
for some $C$ independent of $K$ and that depend only on some norm of $\bbu_{\mathrm h}$, for example the $L^\infty$ norm.

In the particular case of the approximation we consider, we know that
\begin{subequations}
\begin{equation}
\label{appendix:simpson:1}
\frac{1}{\vert K\vert}\int_K \bbu_{\mathrm h}\;\mathrm d\mbx=\xbar\bbu_K=
\sum_{j=1}^6 \beta_j\bbu_{\mathrm h}(\sigma_j)+\beta_7 \bbu_{\mathrm h}(\bbx_K)
\end{equation}
where $\beta_j=\tfrac{9}{20}$ for $j=1,2,3$, $\beta_j=\tfrac{1}{20}$ for $j=4,5,6$ and $\beta_7=\tfrac{2}{15}>0$. We also have
$$\sum_{j=1}^7 \beta_j=1.$$ This is  a Simpson-like formula from which we get
\begin{equation}\label{appendix:simpson:2}
\bbu_{\mathrm h}(\bbx_K)=\theta_7\xbar \bbu_K +\sum_{j=1}^6\theta_j \bbu_{\mathrm h}(\sigma_j)
\end{equation}
\end{subequations}
with $\theta_j=-\tfrac{1}{9}$ for $j=1,2,3$, $\theta_j=-\tfrac{8}{27}$, $j=4,5,6$ and,  $\theta_7=\tfrac{20}{9}$. We also have $$\sum_{j=1}^7 \theta_j=1.$$

Taking $\varphi\in C^1(\R^d)$, we can write
$$
\int_K\varphi \bbv\; \mathrm d\bbx=\vert K\vert \bigg (
 \sum_{j=1}^6\beta_j\varphi(\sigma_j) \bbv(\sigma_j)+\beta_7 \varphi(\bbx_K)\bbv(\bbx_K)\bigg )+O(h^3).
$$
In the following, we set
$$\II(\varphi, \bbv, K)=
 \sum_{j=1}^6\beta_j\varphi(\sigma_j) \bbv(\sigma_j)+\beta_7\varphi(\bbx_K)\bbv(\bbx_K).$$
Using \eqref{appendix:simpson:1} we see that
$$
\II(\varphi, \bbv, K)=\varphi(\bbx_K) \xbar \bbv_K+ \sum_{j=1}^6\beta_j\big ( \varphi(\sigma_j)-\varphi(\bbx_K)\big )\bbv(\sigma_j).$$

We use this on the scheme and get for each $K$
$$\II(\varphi, \dfrac{\mathrm d\bbu_{\mathrm h}}{\mathrm dt}, K)=\varphi(\bbx_K) \dfrac{\mathrm d\xbar\bbu_K}{\mathrm dt}+ \sum_{j=1}^6\beta_j\big ( \varphi(\sigma_j)-\varphi(\bbx_K)\big )\dfrac{\mathrm d\bbu_{\mathrm h}}{\mathrm dt}(\sigma_j),$$
and then
\begin{equation*}
    \begin{aligned}
\II(\varphi, \dfrac{\mathrm d\bbu_{\mathrm h}}{\mathrm dt}, K)&=
-\varphi(\bbx_K) \oint_{\partial K}\bbf(\bbu_{\mathrm h})\cdot \bbn\; \mathrm d\ell\\
&-\sum_{j=1}^6\beta_j\big ( \varphi(\sigma_j)-\varphi(\bbx_K)\big )\bigg (\sum_{K, \sigma_j\in K}\omega_{K,\sigma_j}\Phi_{K,\sigma_j}(\bbu_{\mathrm h})\bigg ).
    \end{aligned}
\end{equation*}
So that, using \eqref{appendix:scheme},  and since \eqref{appendix:scheme:3} also write as
$$\sum_{K, \sigma\in K}\omega_{K,\sigma}\bigg ( \dfrac{\mathrm d\bbu_\sigma}{\mathrm dt}+\Phi_{K,\sigma}(\bbu_{\mathrm h})\bigg )=0$$ by using $\sum_{K, \sigma\in K}\omega_{K,\sigma}=\Id$, we get
\begin{equation}\label{appendix:schema:bis}
\begin{split}
\sum_K  \varphi(\bbx_K)\bigg (& \vert K\vert\dfrac{\mathrm d\xbar\bbu_K}{\mathrm dt}+ \oint_{\partial K}\bbf(\bbu_{\mathrm h})\cdot \bbn\; \mathrm d\ell\bigg )\\+
&\sum_K \sum_{j=1}^6 \beta_j\big ( \varphi(\sigma_j)-\varphi(\bbx_K)\big )\omega_{K,\sigma_j}\bigg (\dfrac{\mathrm d\bbu_{\mathrm h}}{\mathrm dt}(\sigma_j)+
\Phi_{K,\sigma_j}(\bbu_{\mathrm h})\bigg)=0
\end{split}
\end{equation}
The relation \eqref{appendix:schema:bis} is a \emph{rewriting} of the original scheme, tested on any test function $\varphi$.

The (weak) truncation error of the scheme is
\begin{equation}
\label{appendix:weak:TE}
\begin{split}
\varepsilon(\varphi, \bbu_{\mathrm{ex}}):=
&\sum_K \varphi(\bbx_K)\bigg (\vert K\vert  \dfrac{\mathrm d\xbar{(\pi \bbu_{\mathrm{ex}})}}{\mathrm dt}+\oint_{\partial K}\bbf(\bbu_{\mathrm{ex}})\cdot \bbn\;\mathrm d\ell\bigg ) \\
&+
\sum_K\sum_{j=1}^6 {\beta_j} \vert K\vert  \big ( \varphi(\sigma_j)-\varphi(\bbx_K)\big ){\omega_{K,\sigma_j}}\bigg ( \dfrac{\mathrm d\bbu_{\mathrm{ex}}}{\mathrm dt}(\sigma_j)+\Phi_{K,\sigma_j}(\pi\bbu_{\mathrm{ex}})\bigg )
\end{split}
\end{equation}
where 
$$\pi\bbu_{\mathrm{ex}}=\sum_{j=1}^6 \bbu_{\mathrm{ex}}(\sigma_i)\varphi_{\sigma_i}+\xbar \bbu_{\mathrm{ex}}\bar \varphi.$$

Now by taking $\Phi_{K,\sigma_j}$ as in \eqref{appendix:zeschema}, 
we first have that
$${\vert K\vert} \PP^{-1}\dfrac{\mathrm d}{\mathrm dt}\begin{pmatrix} \vdots \\ \bbu_{\mathrm{ex}}(\sigma_j)\\  \vdots \\ \xbar\bbu_{\mathrm{ex}}\end{pmatrix}+
\begin{pmatrix} \vdots \\ -\int_{K}\nabla\varphi_{\sigma_l}\cdot\bbf(\bbu_{\mathrm{ex}})\;\mathrm d\bbx+\oint_{\partial K}\varphi{\sigma_l}\bbf(\bbu_{\mathrm{ex}})\cdot \bbn \; \mathrm d\ell \\ \vdots \\-\int_{K}\nabla\xbar \varphi\cdot\bbf(\bbu_{\mathrm{ex}})+\oint_{\partial K}\xbar \varphi\bbf(\bbu_{\mathrm{ex}})\cdot \bbn \; \mathrm d\ell\end{pmatrix}=0,$$
so that we can write, writing $e=\pi\bbu_{\mathrm{ex}}-\bbu_{\mathrm{ex}}$,
\begin{equation*}
\begin{split}\varepsilon(\varphi, \bbu_{\mathrm{ex}})&=
\sum_K\varphi(\bbx_K) \bigg (\vert K\vert  \dfrac{\mathrm d\xbar e_K}{\mathrm dt}+\oint_{\partial K}\big (\bbf(\bbu_{\mathrm h})-\bbf(\pi(\bbu_{\mathrm {ex}})\big )\cdot \bbn\;\mathrm d\ell\big ) \\
&\hspace{-1.5em}+
\sum_K\sum_{j=1}^6 {\beta_j\omega_{K,\sigma_j}} \vert K\vert  \big ( \varphi(\sigma_j)-\varphi(\bbx_K)\big )\bigg (\underbrace{ \dfrac{\mathrm de_{\sigma_j}}{\mathrm dt}+\big (\Phi_{K,\sigma_j}(\pi\bbu)-\Phi_{K,\sigma_j}(\pi\bbu_{\mathrm {ex}})\big )}_{(II)}\bigg ).
\end{split}
\end{equation*}
Since $\Vert e\Vert_{\infty}\leq C\; h^3$ where the constant $C$ only depends on the $L^\infty$ norm of $\bbu_{\mathrm {ex}}$ and its derivatives, we see that 
$\tfrac{\mathrm d\xbar e_K}{\mathrm dt}$ is $O(h^3)$, while 
\begin{equation*}
\begin{split}\sum_K \varphi(\bbx_K)&\oint_{\partial K}\big (\bbf(\bbu_{\mathrm h})-\bbf(\pi(\bbu_{\mathrm {ex}})\big )\cdot \bbn\;\mathrm d\ell\\
&=\sum_{\text{edges} }\big (\varphi(\bbx_K)-\varphi(\bbx_{K'})\big )\int_e\big (\bbf(\bbu_{\mathrm h})-\bbf(\pi(\bbu_{\mathrm {ex}})\big )\cdot \bbn\;\mathrm d\ell=\vert K\vert O(h^{3})
\end{split}
\end{equation*}
since $\varphi(\bbx_K)-\varphi(\bbx_{K'})=O(h)$, so that $\vert e\vert (\varphi(\bbx_K)-\varphi(\bbx_{K'}) =O(h^2)=\vert K\vert O(1)$ for a regular mesh. This leads to  
$$\sum_K\varphi(\bbx_K) \bigg (\vert K\vert  \dfrac{\mathrm d\xbar e_K}{\mathrm dt}+\oint_{\partial K}\big (\bbf(\bbu_{\mathrm h})-\bbf(\pi(\bbu_{\mathrm {ex}})\big )\cdot \bbn\;\mathrm d\ell\bigg ) =O(h^3).$$
The term $(II)$ is only $O(h^2)$ because $\Phi_{K,\sigma_j}(\pi\bbu)-\Phi_{K,\sigma_j}(\pi\bbu_{\mathrm {ex}})=O(h^2)$ but since $\varphi(\sigma_j)-\varphi(\bbx_K)=O(h)$, we see that in the end
$$\sum_K\sum_{j=1}^6 \beta_j\omega_{K,\sigma_j} \vert K\vert  \big ( \varphi(\sigma_j)-\varphi(\bbx_K)\big )\bigg ( \dfrac{\mathrm de_{\sigma_j}}{\mathrm dt}+\big (\Phi_{K,\sigma_j}(\pi\bbu)-\Phi_{K,\sigma_j}(\pi\bbu_{\mathrm{ex}})\big )\bigg )=O(h^3),$$
so that $\varepsilon(\varphi, \bbu_{\mathrm{ex}})=O(h^3)$ for regular meshes.

\begin{remark} This is not specific to two dimensional problems, because we have used that $\vert K\vert =O(h^2)$ and $\vert e\vert=O(h)$, that is in general $\vert K\vert =O(h^d)$ and $\vert e\vert =O(h^{d-1})$ for a regular mesh.
\end{remark}

\end{document}